\documentclass[leqno,12pt]{amsart}
\usepackage{amsmath,amstext,amssymb,amsopn,amsthm,mathrsfs}

\theoremstyle{plain}
\newtheorem{teo}{Theorem}[section]
\newtheorem{defi}{Definition}[section]

\newtheorem{lemma}{Lemma}[section]
\newtheorem{proposition}{Proposition}[section]
\numberwithin{equation}{section}
\newcommand{\dem}{\medskip \par \noindent \mbox{\bf Proof. }}
\def\ep{\hfill{$\Box $}}

\begin{document}

\title[Riesz \& Bessel Potentials and Fractional Derivatives ] {Riesz Potentials, Bessel Potentials and Fractional Derivatives  on Besov-Lipschitz spaces for the Gaussian Measure.}

\author{A. Eduardo Gatto}
\address{Department of Mathematical Sciences, DePaul University, Chicago, IL,
   60614, USA.}
   \email{aegatto@depaul.edu}
   \author{Ebner Pineda}
\address{Departamento de Matem\'{a}tica,  Decanato de Ciencia y
Tecnologia, UCLA
 Apartado 400 Barquisimeto 3001 Venezuela.}
\email{epineda@uicm.ucla.edu.ve}
\author{Wilfredo~O.~Urbina}
\address{Department of Mathematics and Actuarial Sciences, Roosevelt University, Chicago, IL,
   60605, USA.}
\email{wurbinaromero@roosevelt.edu}



\begin{abstract}
In \cite{gatur} Gaussian Lipschitz spaces $Lip_{\alpha}(\gamma_d)$ were considered and then the  boundedness properties of Riesz Potentials, Bessel potentials and Fractional
Derivatives were studied in detail. In this paper we will study the boundedness of those operators on Gaussian Besov-Lipschitz spaces $B_{p,q}^{\alpha}(\gamma_d)$. Also these results can be extended to the case of Laguerre or Jacobi expansions and even further to the general framework of diffusions semigroups.
\end{abstract}

\maketitle

\section{Introduction}

On $\mathbb{R}^d$ let us consider the Gaussian measure 
\begin{equation}
\gamma_d(x)=\frac{e^{-\left|x\right|^2}%
}{\pi^{d/2}} dx, \, x\in\mathbb{R}^d
\end{equation}
 and the Ornstein-Uhlenbeck
differential operator
\begin{equation}\label{OUop}
L=\frac12\triangle_x-\left\langle x,\nabla _x\right\rangle.
\end{equation}

Let $\nu=(\nu _1,...,\nu_d)$ be a
multi-index such that $\nu _i \geq 0, i= 1, \cdots, d$,  let $\nu
!=\prod_{i=1}^d\nu _i!,$ $\left| \nu \right| =\sum_{i=1}^d\nu _i,$ $%
\partial _i=\frac \partial {\partial x_i},$ for each $1\leq i\leq d$ and $%
\partial ^\nu =\partial _1^{\nu _1}...\partial _d^{\nu _d},$ consider the normalized Hermite polynomials of order $\nu$ in $d$ variables,
\begin{equation}
h_\nu (x)=\frac 1{\left( 2^{\left| \beta \right| }\nu
!\right)
^{1/2}}\prod_{i=1}^d(-1)^{\nu _i}e^{x_i^2}\frac{\partial ^{\nu _i}}{%
\partial x_i^{\nu _i}}(e^{-x_i^2}),
\end{equation}
it is well known, that the Hermite polynomials are
eigenfunctions of the operator $L$,
\begin{equation}\label{eigen}
L h_{\nu}(x)=-\left|\nu \right|h_\nu(x).
\end{equation}
Given a function $f$ $\in L^1(\gamma _d)$ its
$\nu$-Fourier-Hermite coefficient is defined by
\[
\hat{f}(\nu) =<f, h_\nu>_{\gamma_d}
=\int_{\mathbb{R}^d}f(x)h_\nu (x)\gamma _d(dx).
\]
Let $C_n$ be the closed subspace of $L^2(\gamma_d)$ generated by
the linear combinations of $\left\{ h_\nu \ :\left| \nu
\right| =n\right\}$. By the orthogonality of the Hermite
polynomials with respect to $\gamma_d$ it is easy to see that
$\{C_n\}$ is an orthogonal decomposition of $L^2(\gamma_d)$,
$$ L^2(\gamma_d) = \bigoplus_{n=0}^{\infty} C_n,$$
this decomposition is called the Wiener chaos.

Let $J_n$ be the orthogonal projection  of $L^2(\gamma_d)$ onto
$C_n$, then if $f\in L^2(\gamma_d)$
\[
J_n f=\sum_{\left|\beta\right|=n}\hat{f}(\nu) h_\nu.
\]
Let us define the Ornstein-Uhlenbeck semigroup $\left\{
T_t\right\} _{t\geq 0}$ as
\begin{eqnarray}\label{01}
\nonumber T_t f(x)&=&\frac 1{\left( 1-e^{-2t}\right) ^{d/2}}\int_{\mathbb{R}^d}e^{-\frac{%
e^{-2t}(\left| x\right| ^2+\left| y\right| ^2)-2e^{-t}\left\langle
x,y\right\rangle }{1-e^{-2t}}}f(y)\gamma _d(dy)\\
& = & \frac{1}{\pi^{d/2}(1-e^{-2t})^{d/2}}\int_{\mathbb R^d} e^{-
\frac{|y-e^{-t}x|^2}{1-e^{-2t}}} f(y) dy
\end{eqnarray}
The family $\left\{ T_t\right\}_{t\geq 0}$ is a strongly
continuous Markov semigroup on $L^p(\gamma_d)$, $1 \leq p <
\infty$, with infinitesimal generator $L$. Also, by a change of
variable we can write,
\begin{equation}\label{t1}
T_t f(x)=\int_{\mathbb{R}^d} f(\sqrt{1-e^{-2t}}u + e^{-t}x)\gamma
_d(du).
\end{equation}

Now, by Bochner subordination formula, see Stein \cite{se70} page 61, we
define the Poisson-Hermite semigroup $\left\{ P_t\right\} _{t\geq
0}$ as
\begin{equation}\label{PoissonH}
 P_t f(x)=\frac 1{\sqrt{\pi }}\int_0^{\infty} \frac{e^{-u}}{\sqrt{u}}T_{t^2/4u}f(x)du
\end{equation}
From (\ref{01}) we
obtain, after the change of variable $r=e^{-t^2/4u}$,
\begin{eqnarray}\label{03}
\nonumber P_t f(x)&=&\frac 1{2\pi
^{(d+1)/2}}\int_{\mathbb{R}^d}\int_0^1t\frac{\exp \left( t^2/4\log
r\right) }{(-\log r)^{3/2}}\frac{\exp \left( \frac{-\left|
y-rx\right| ^2}{1-r^2}\right) }{(1-r^2)^{d/2}}\frac{dr}rf(y)dy\\
&=& \int_{\mathbb{R}^d} p(t,x,y) f(y)dy,
\end{eqnarray}
with
\begin{equation}
p(t,x,y) = \frac 1{2\pi ^{(d+1)/2}}\int_0^1t\frac{\exp \left(
t^2/4\log r\right) }{(-\log r)^{3/2}}\frac{\exp \left(
\frac{-\left| y-rx\right| ^2}{1-r^2}\right)
}{(1-r^2)^{d/2}}\frac{dr}r.
\end{equation}
Also by the change of variable $s= t^2/4u$ we have,
\begin{equation}
P_t f(x)=\frac 1{\sqrt{\pi }}\int_0^{\infty} \frac{e^{-u}}{\sqrt{u}}T_{t^2/4u}f(x)du
=\int_0^{\infty} T_s f(x) \mu^{(1/2)}_t(ds),
\end{equation}
where the measure
\begin{equation}\label{onesided1/2}
\mu^{(1/2)}_t(ds) = \frac t{2\sqrt{\pi
}}\frac{e^{-t^2/4s}}{s^{3/2}}ds,
\end{equation}
is called the one-side stable measure on $(0, \infty)$ of order
$1/2$.

The family $\left\{ P_t\right\}_{t\geq 0}$ is also a strongly
continuous semigroup on $L^p(\gamma_d)$, $1 \leq p < \infty$,
with infinitesimal generator $-(-L)^{1/2}$. In what follows, often we are  going to use the notation $$u(x,t) = P_{t}f(x),$$
and
$$ u^{(k)}(x,t) = \frac{\partial^k}{\partial t^k} P_{t}f(x).$$

Observe that by (\ref{eigen}) we have that
\begin{equation}\label{OUHerm}
T_t h_\nu(x)=e^{-t\left| \nu\right|}h_\nu(x),
\end{equation}
and
\begin{equation}\label{PHHerm}
 P_t h_\nu(x)=e^{-t\sqrt{\left| \nu\right|}}h_\nu(x),
\end{equation}
i.e. the Hermite polynomials are eigenfunctions of $T_t$ and $P_t$ for any $t \geq 0$.\\

The operators that we are going to consider in this paper are the following:
\begin{itemize}

\item  For $\beta>0$, the Fractional Integral or Riesz potential of order $\beta$, $I_\beta^{\gamma}$, with respect to the Gaussian measure is defined formally as
\begin{equation}\label{i1}
I_\beta=(-L)^{-\beta/2}\Pi_{0},
\end{equation}
where,
$\Pi_{0}f=f-\displaystyle\int_{\mathbb{R}^{d}}f(y)\gamma_{d}(dy)$,
for $f\in L^{2}(\gamma_{d})$. That means that for  the Hermite
polynomials $\{h_\beta\}$, for $\left|\beta\right|>0$,
\begin{equation}\label{e4}
I_\beta h_\nu(x)=\frac 1{\left|
\nu \right|^{\beta/2}}h_\nu(x),
\end{equation}
and for $\nu=\overline{0}, \,
I_{\beta}(h_{\overline{0}})=0.$ Then by linearity can be extended to any polynomial.

Now, it is easy to see that if $f$ is a polynomial,
\begin{equation}\label{e3}
I_\beta f(x)  =\frac 1{\Gamma(\beta)}\int_0^{\infty}
t^{\beta-1}(P_t f(x) -P_{\infty} f(x))\,dt.
\end{equation}
Moreover  by P. A. Meyer's multiplier theorem, see \cite{me3} or \cite{wat}, $I_\alpha$ admits a continuous extension  to $ L^p(\gamma _d)$,  $1 < p < \infty$, and then (\ref{e3}) can be extended for $f \in L^{p}(\gamma_d)$.\\

\item  The Bessel Potential of order $\beta>0,$
$\mathcal{J}_\beta$, associated to the Gaussian  measure
is defined formally as
\begin{eqnarray}
\mathcal{J}_\beta= (I+\sqrt{-L})^{-\beta},
\end{eqnarray}
meaning that for the Hermite polynomials we have,
\begin{eqnarray*}
\mathcal{J}_\beta h_\nu(x)=\frac 1{(1+\sqrt{\left|
\nu\right|})^{\beta}}h_\nu(x).
\end{eqnarray*}
Again  by linearity can be extended to any polynomial and Meyer's theorem allows us to extend Bessel Potentials
to a \mbox{continuous} operator on $L^p(\gamma_d),$ $1 < p <
\infty$.  It can be proved that the Bessel potentials can be
represented as
\begin{equation}\label{Beselrepre}
\mathcal{J}_\beta
f(x)=\frac{1}{\Gamma(\beta)}\int_{0}^{+\infty}t^{\beta}e^{-t}P_{t}f(x)\frac{dt}{t}.
\end{equation}
Moreover $\{\mathcal{J}_\beta\}_\beta$
is a strongly continuous semigroup on $L^p(\gamma_d)$, $1 \leq p <\infty$, with infinitesimal generator $\frac{1}{2}\log (I-L),$ see \cite{forscotur}. \\

\item The Riesz fractional derivate of order $\alpha >0$ with respect to the
Gaussian measure $D^\alpha$, is defined  formally as
\begin{equation}
D^\beta=(-L)^{\beta/2},
\end{equation}
meaning that for the Hermite polynomials, we have
\begin{equation}\label{e6}
D^\beta(x)=\left| \nu\right|^{\beta/2}
h_\nu(x),
\end{equation}
 thus by linearity can be extended to any polynomial.

  The Riesz fractional derivate $D^\beta$  with respect to the Gaussian
measure was first introduced in \cite{lu}. For more detail we
refer to that article. Also see \cite{ebner} for improved  and
simpler proofs of some results contained there.
In the case of $0 < \beta <1$ we have the following integral representation,  
\begin{equation}\label{e5}
D^\beta f =\frac 1{c_\beta}\int_0^{\infty}t^{-\beta-1}(P_t  -I)\, f dt,
\end{equation} 
where
$
c_\beta=\int_0^\infty u^{-\beta-1}(e^{-u}-1)du.
$ Moreover for  $f$ $\in C_B^2(\mathbb{R}^d)$, i.e. the set of  two times continuously differentiable functions with bounded derivatives,  then it can be proved using integration by parts (for details see \cite{lu}), that 
\begin{equation} \label{p1}
D^{\beta} f =\frac {1}{\beta c_{\beta}}\int_0^{\infty} t^{-\beta}\frac{\partial}{\partial t}P_t f dt. \end{equation}

Moreover, if $\beta\geq1$, let  $k$   be the smallest integer greater than $\beta$ i.e. $ k-1 \leq \beta < k$, then the fractional derivative $D^\beta$ can be represented as 
\begin{equation}
D^\beta f = \frac{1}{c^k_{\beta}}\int_0^{\infty} t^{-\beta-1} ( P_t -I )^k f \, dt,
\end{equation}
where $c^k_{\beta} = \int_0^{\infty} u^{-\beta-1} ( e^{-u}-1 )^k  \, du$. Now, if $f$ is a polynomial, by the linearity of the operators
$I_\beta$ and $D^\beta$, (\ref{e4}) and
(\ref{e6}), we get
\begin{equation}\label{di}
 \Pi_0f=I_\beta(D^\beta f)=D^\beta(I_\beta f).\\
\end{equation}

\item We can also consider a Bessel fractional derivative ${\mathcal D}^\beta$,  defined formally as
\[
{\mathcal D}^\beta=(I+\sqrt{-L})^{\beta},
\]
which means that for the Hermite polynomials, we have
\begin{equation}\label{e7}
{\mathcal D}^\beta h_\nu(x)=(1+ \sqrt{\left|
\nu\right|}))^{\beta} h_\nu(x),
\end{equation}
In the case of $0 < \beta <1$ we have the following integral representation,  
\begin{equation}\label{e8}
{\mathcal D}^\beta f =\frac 1{c_\beta}\int_0^{\infty}t^{-\beta-1}( e^{-t} P_t -I) \, f dt,
\end{equation} 
where, as before, 
$
c_\beta=\int_0^\infty u^{-\beta-1}(e^{-u}-1)du.$
Moreover,  if $\beta >1$ let $k$  be the smallest integer greater than $\beta$ i.e. $ k-1 \leq \beta < k$, then the fractional derivative ${\mathcal D}^\beta$ can be represented  as 
\begin{equation}
{\mathcal D}^\beta f =  \frac{1}{c^k_{\beta}} \int_0^{\infty} t^{-\beta-1} (e^{-t} P_t -I )^k\,  f \, dt,
\end{equation}
where $c^k_{\beta} =  \int_0^{\infty} u^{-\beta-1} (e^{-u} -1 )^k du.$\\

\end{itemize}

The Gaussian Besov-Lipschitz $B_{p,q}^{\alpha}(\gamma_d)$
spaces were introduced in \cite{piur02}, see also \cite{ebner}, as follows

\begin{defi}
Let $\alpha > 0$, $k$ be the smallest integer greater than
$\alpha$, and $1\leq p,q \leq\infty$. For $1 \leq q < \infty$ the Gaussian Besov-Lipschitz
space $B_{p,q}^{\alpha}(\gamma_d)$ are defined as the set of functions  $f \in
L^p(\gamma_d)$ for which
\begin{equation}\label{e15}
 \left( \int_0^{\infty} (t^{k-\alpha} \left\|
\frac{\partial^{k}P_t f}{\partial t^{k}} \right\|_{p,\gamma_d}
)^{q} \frac{dt}{t} \right) ^{1/q}  < \infty.
\end{equation}
The norm of $f \in B_{p,q}^{\alpha}(\gamma_d)$ is defined as
\begin{equation}
\left\| f \right\|_{B_{p,q}^{\alpha}}: =  \left\| f \right\|_{p,
\gamma_d} +\left( \int_0^{\infty} (t^{k-\alpha} \left\|
\frac{\partial^{k}   P_t f}{\partial t^{k}} \right\|_{p,\gamma_d}
)^{q} \frac{dt}{t} \right) ^{1/q}
\end{equation}
 For $q=\infty$ the Gaussian Besov-Lipschitz space
$B_{p,\infty}^{\alpha}(\gamma_d)$ are defined  as the set of  functions  $f \in
L^p(\gamma_d)$ for which exists a constant $A$ such that\\
$$\|\frac{\partial^{k}P_t f}{\partial
t^{k}}\|_{p,\gamma_d}\leq At^{-k+\alpha}$$ and then the norm of $f
\in B_{p,\infty}^{\alpha}(\gamma_d)$ is defined as
\begin{equation}
\left\| f \right\|_{B_{p,\infty}^{\alpha}}: =  \left\| f
\right\|_{p, \gamma_d} +A_{k}(f),
\end{equation}
where $A_{k}(f)$ is the smallest constant $A$ appearing in the
above inequality. In particular, the space $B_{\infty,\infty}^{\alpha}(\gamma_d)$ is the Gaussian Lipschitz space $Lip_{\alpha}(\gamma_d)$.
\end{defi}

The definition of $B_{p,q}^{\alpha}(\gamma_d)$ does not depend on which $k>\alpha$ is chosen and the resulting norms are equivalent,  for the proof of this result and other properties of these spaces see \cite{piur02}.\\

In what follows, we need the following
technical result about $L^p(\gamma_d)$-norms of the derivatives of
the Poisson-Hermite semigroup, see \cite{piur02}, Lemma 2.2

\begin{lemma}\label{kdecay}
Suppose $f\in L^{p}(\gamma_{d})$, $1 \leq p < \infty$ then for any integer $k$ the function
$\displaystyle\|\frac{\partial^{k}P_t f}{\partial
t^{k}}\|_{p,\gamma_{d}}$ is a non-increasing function of $t$, for
$0<t<+\infty$. Moreover,
\begin{equation}\label{kdecayine}
\|\frac{\partial^{k}P_t f}{\partial t^{k}}\|_{p,\gamma_d}\leq
C \| f\|_{p,\gamma_d}t^{-k}, t>0
\end{equation}

\end{lemma}

Also we will need some inclusion relations among the Gaussian Besov-Lipschitz
spaces, see \cite{piur02}, 

\begin{proposition}\label{incluBesov}
The inclusion $B_{p,q_1}^{\alpha_{1}}(\gamma_d)\subset
B_{p,q_2}^{\alpha_{2}}(\gamma_d)$ holds if either:
\begin{enumerate}
\item [i)]   $\alpha_{1}>\alpha_{2}>0$ where $q_{1}$ and $q_{2}$ need
not to be related, or
 \item[ii)] If $\alpha_{1}=\alpha_{2}$ and $q_{1}\leq q_{2}$.
\end{enumerate}
\end{proposition}
In \cite{gatur} Gaussian Lipschitz spaces $Lip_{\alpha}(\gamma_d)$ were considered and the boundedness of  Riesz Potentials, Bessel potentials and Fractional
Derivatives on them were study. In the next section, we are going to extend those results for Gaussian Besov-Lipschitz spaces, but not including them. Thus, the main purpose of this paper is to study the boundedness of Gaussian fractional integrals and derivatives associated to Hermite polynomial expansions on  Gaussian Besov-Lipschitz spaces  $B_{p,q}^{\alpha}(\gamma_{d})$ . To get these results we introduce formulas for these operators in terms of the Hermite-Poisson semigroup as well as the Gaussian Besov-Lipschitz spaces. This approach was originally developed for the classical  Poisson integral, see Stein  \cite{se70}, Chapter V Section 5. These proofs  can also be extended to the case of Laguerre and Jacobi expansions. These results can be also obtained using abstract interpolation theory on the the Poisson-Hermite semigroup, see \cite{trie2}.\\

As usual in what follows $C$ represents a constant that is not necessarily the same in each occurrence.


\section{Main results}

 In the case of the Lipschitz spaces only a truncated version of the Riesz potentials is bounded from 
$Lip_{\alpha}(\gamma_d)$ to  $Lip_{\alpha+\beta}(\gamma_d)$, see \cite{gatur} Theorem 3.2. Now, we wil study the boundedness properties of the Riesz potentials on  Besov-Lipschitz spaces, and we will see that in this case the results are actually better.

 \begin{teo}
Let $\alpha\geq 0, \beta>0$, $1<p<\infty, 1\leq q \leq \infty$ then $I_{\beta}$ is
bounded from $B_{p,q}^{\alpha}(\gamma_{d})$ into
$B_{p,q}^{\alpha+\beta}(\gamma_{d})$.
\end{teo}
\dem

 Let $k>\alpha+\beta$ a fixed integer,  $f\in B_{p,q}^{\alpha}(\gamma_{d})$,  using the integral representation of Riesz Potentials (\ref{e3}), the semigroup property of $\{P_{t}\}$ and the fact that $P_{\infty}f(x)$ is a constant and the semigroup is conservative, we get
\begin{eqnarray}\label{est1}
\nonumber P_{t}I_{\beta}f(x)&=&\displaystyle\frac{1}{\Gamma(\beta)}\int_{0}^{+\infty}s^{\beta-1}P_{t}(P_{s}f(x)-P_{\infty}f(x))ds\\
&=&\displaystyle\frac{1}{\Gamma(\beta)}\int_{0}^{+\infty}s^{\beta-1}(P_{t+s}f(x)-P_{\infty}f(x))ds.
\end{eqnarray}
Now using again that $P_{\infty}f(x)$ is a constant and the chain rule,
\begin{eqnarray}\label{est2}
\nonumber \frac{\partial^{k}}{\partial t^{k}}(P_{t}I_{\beta}f)(x) &=&\frac{1}{\Gamma(\beta)}\int_{0}^{+\infty}s^{\beta-1}\frac{\partial^{k} }{\partial t^{k}}(P_{t+s}f(x)-P_{\infty}f(x))\, ds\\
&=& \frac{1}{\Gamma(\beta)}\int_{0}^{+\infty}s^{\beta-1}u^{(k)}(x,t+s) \,ds.
\end{eqnarray}
Now, by Minkowski's integral inequality 
\begin{eqnarray}\label{est3}
\|\frac{\partial^{k}}{\partial
t^{k}} P_{t}I_{\beta}f \|_{p,\gamma}&\leq&\frac{1}{\Gamma(\beta)}\int_{0}^{+\infty}s^{\beta-1}\|u^{(k)}(\cdot,t+s)\|_{p,\gamma}ds  .
\end{eqnarray}
Then, if $ 1 \leq q < \infty$,
\begin{eqnarray*}
&&\displaystyle\big(\int_{0}^{+\infty}\big(t^{k-(\alpha+\beta)}\|\frac{\partial^{k}
}{\partial
t^{k}}(P_{t}I_{\beta}f) \|_{p,\gamma}\big)^{q}\frac{dt}{t}\big)^{\frac{1}{q}}
\quad \quad\quad \quad
\quad \quad\quad\\
&&\quad \quad \quad \leq \frac{1}{\Gamma(\beta)}\displaystyle\big(\int_{0}^{+\infty}t^{(k-(\alpha+\beta))q}\big(\int_{0}^{+\infty}s^{\beta-1}\|u^{(k)}(\cdot,t+s)\|_{p,\gamma}ds\big)^{q}\frac{dt}{t}\big)^{\frac{1}{q}}\\
&& \quad \quad \quad  \leq C_{\beta}\displaystyle\big(\int_{0}^{+\infty}t^{(k-(\alpha+\beta))q}\big(\int_{0}^{t}s^{\beta-1}\|u^{(k)}(\cdot,t+s)\|_{p,\gamma}ds\big)^{q}\frac{dt}{t}\big)^{\frac{1}{q}}\\
&& \quad \quad \quad \quad \quad  + C_{\beta}\displaystyle\big(\int_{0}^{+\infty}t^{(k-(\alpha+\beta))q}\big(\int_{t}^{+\infty}s^{\beta-1}\|u^{(k)}(\cdot,t+s)\|_{p,\gamma}ds\big)^{q}\frac{dt}{t}\big)^{\frac{1}{q}}\\
&&\quad \quad \quad  = (I)+(II).
\end{eqnarray*}
Now, as $\beta>0$ using Lemma \ref{kdecay}, as  $t+s>t$,
\begin{eqnarray*}
(I) &\leq& C_{\beta} \big(\int_{0}^{+\infty}t^{(k-(\alpha+\beta))q}\big(\int_{0}^{t}s^{\beta-1}
\|u^{(k)}(\cdot,t)\|_{p,\gamma}ds\big)^{q}\frac{dt}{t}\big)^{\frac{1}{q}}\\
&=&C_{\beta} \big(\int_{0}^{+\infty}t^{(k-(\alpha+\beta))q}\|\frac{\partial^{k}P_{t}f}{\partial t^{k}}\|_{p,\gamma}^{q}(\frac{t^{\beta}}{\beta})^{q}\frac{dt}{t}\big)^{\frac{1}{q}}\\
&=&C'_{\beta} \big(\int_{0}^{+\infty}\big(t^{k-\alpha}\|\frac{\partial^{k}P_{t}f}{\partial t^{k}}\|_{p,\gamma}\big)^{q}\frac{dt}{t}\big)^{\frac{1}{q}}<+\infty,
\end{eqnarray*}
since $f\in B_{p}^{\alpha,q}(\gamma_{d})$. \\

On the other hand, as $k>\alpha+\beta$ using again Lemma  \ref{kdecay}, since $t+s>s$, and Hardy's inequality (\ref{hardy2}),  we obtain
\begin{eqnarray*}
(II) &\leq&C_{\beta}\big(\int_{0}^{+\infty}t^{(k-(\alpha+\beta))q}\big(\int_{t}^{+\infty}s^{\beta}
\|u^{(k)} (\cdot, s)\|_{p,\gamma}\frac{ds}{s}\big)^{q}\frac{dt}{t}\big)^{\frac{1}{q}}\\
&\leq&\frac{C_{\beta}}{k-(\alpha+\beta)}\int_{0}^{+\infty}\big(s^{k-\alpha}
\|\frac{\partial^{k} P_{s}f }{\partial s^{k}}\|_{p,\gamma}\big)^{q}\frac{ds}{s}\big)^{\frac{1}{q}}<+\infty,
\end{eqnarray*}
since $f\in B_{p,q}^{\alpha}(\gamma_{d})$. Therefore $I_{\beta}f\in
B_{p,q}^{\alpha+\beta}(\gamma_{d})$ and moreover,
\begin{eqnarray*}
\displaystyle\|I_{\beta}f\|_{B_{p,q}^{\alpha+\beta}}&=&\|I_{\beta}f\|_{p,\gamma}+\displaystyle\big(\int_{0}^{+\infty}\big(t^{k-(\alpha+\beta)}\|\frac{\partial^{k}
}{\partial
t^{k}}(P_{t}I_{\beta}f)\|_{p,\gamma}\big)^{q}\frac{dt}{t}\big)^{\frac{1}{q}}\\
&\leq&C\|f\|_{p,\gamma}+C _{\alpha,\beta}\big(\int_{0}^{+\infty}\big(t^{k-\alpha}\|\frac{\partial^{k}P_{t}f}{\partial t^{k}}\|_{p,\gamma}\big)^{q}\frac{dt}{t}\big)^{\frac{1}{q}}\\
&\leq&C\displaystyle\|f\|_{B_{p,q}^{\alpha}}.
\end{eqnarray*}

Now if $q=\infty$, (\ref{est3}) can be written as
\begin{eqnarray*}
\|\frac{\partial^{k}}{\partial
t^{k}} P_{t}I_{\beta}f \|_{p,\gamma}&\leq&\frac{1}{\Gamma(\beta)}\int_{0}^{+\infty}s^{\beta-1}\|u^{(k)}(\cdot,t+s)\|_{p,\gamma}ds\\
&=& \frac{1}{\Gamma(\beta)}\int_{0}^{t}s^{\beta-1}\|u^{(k)}(\cdot,t+s)\|_{p,\gamma}ds\\
&& \quad \quad \quad + \frac{1}{\Gamma(\beta)}\int_{t}^{\infty}s^{\beta-1}\|u^{(k)}(\cdot,t+s)\|_{p,\gamma}ds\\
&=&(I)+(II) .
\end{eqnarray*}
Now,  using that $\beta>0$, Lemma \ref{kdecay},  as $t+s>t $ and since
$f\in B_{p,\infty}^{\alpha}(\gamma_{d})$,
\begin{eqnarray*}
(I) &\leq&\frac{1}{\Gamma(\beta)}\|\frac{\partial^{k}P_{t}f}{\partial
t^{k}}\|_{p,\gamma}\int_{0}^{t}s^{\beta-1}ds \leq \frac{1}{\Gamma(\beta)}\frac{t^{\beta}}{\beta}A_{k}(f)t^{-k+\alpha}\\
&=&C_{\beta} A_{k}(f)\, t^{-k+\alpha+\beta} .
\end{eqnarray*}
On the other hand, since $k>\alpha+\beta$, using Lemma
 \ref{kdecay}, as $t+s>s$ and since $f\in B_{p,\infty}^{\alpha}(\gamma_{d})$, we get
\begin{eqnarray*}
(II) &\leq&\frac{1}{\Gamma(\beta)}\int_{t}^{\infty}s^{\beta-1}\|\frac{\partial^{k}P_{s}f}{\partial
s^{k}}\|_{p,\gamma}ds
\leq \frac{A_{k}(f)}{\Gamma(\beta)}\int_{t}^{\infty}s^{-k+\alpha+\beta-1}ds\\
&=&\frac{A_{k}(f)}{\Gamma(\beta)}\frac{t^{-k+\alpha+\beta}}{k-(\alpha+\beta)} = C_{k, \alpha,\beta} \, t^{-k+\alpha+\beta}.
\end{eqnarray*}
Therefore
\begin{eqnarray*}
\|\frac{\partial^{k}}{\partial t^{k}}
P_{t}I_{\beta}f\|_{p,\gamma}&\leq&C A_{k}(f)t^{-k+\alpha+\beta},\quad t>0,
\end{eqnarray*}
and this implies that
$I_{\beta}f \in B_{p,\infty}^{\alpha+\beta}(\gamma_{d})$ and $A_{k}(I_{\beta}f)\leq C A_{k}(f)$.\\
Moreover, as $I_{\beta}$ is bounded operator on $L^{p}(\gamma_{d}), 1<p<\infty$,
\begin{eqnarray*}
\|I_{\beta}f\|_{B_{p,\infty}^{\alpha+\beta}}&=&\|I_{\beta}f\|_{p,\gamma}+A_{k}(I_{\beta}f)\\
&\leq&\|f\|_{p,\gamma}+CA_{k}(f) \leq C \|f\|_{B_{p,\infty}^{\alpha}}.
\end{eqnarray*}
\ep \\


Now we want to study the boundedness properties of the Bessel potentials on Besov-Lipschitz spaces. In \cite{gatur}, Theorem 3.1,  the following result was proved,
\begin{teo}
Let $\alpha\geq 0, \beta>0$ then
 $\mathcal{J}_{\beta}$ is
bounded from $Lip_{\alpha}(\gamma_{d})$ into
$Lip_{\alpha+\beta}(\gamma_{d})$.
\end{teo}

Also in \cite{piur02}, Theorem 2.4,  it was proved that 
\begin{teo}
Let $\alpha\geq 0, \beta>0$ then for $1\leq p, q < \infty$
$\mathcal{J}_{\beta}$ is
bounded from $B_{p,q}^{\alpha}(\gamma_{d})$ into
$B_{p,q}^{\alpha+\beta}(\gamma_{d})$.
\end{teo}

Therefore the following result is the only case that was missing,

\begin{teo}
Let $\alpha\geq 0, \beta>0$ then for $1\leq p < \infty$
 $\mathcal{J}_{\beta}$ is
bounded from $B_{p,\infty}^{\alpha}(\gamma_{d})$ into
$B_{p,\infty}^{\alpha+\beta}(\gamma_{d})$.
\end{teo}
\dem

Let $k>\alpha+\beta$ a fixed integer, $f\in
B_{p,\infty}^{\alpha}(\gamma_{d})$,
by using the representation of Bessel potential (\ref{Beselrepre}), we get 
$$P_{t}(\mathcal{J}_{\beta}f)(x)=\displaystyle\frac{1}{\Gamma(\beta)}\int_{0}^{+\infty}s^{\beta}e^{-s}P_{t+s}f(x)\frac{ds}{s},$$
thus using the chain rule, we obtain
\begin{eqnarray*}
\frac{\partial^{k}}{\partial t^{k}}P_{t}(\mathcal{J}_{\beta}f)(x)&=&\frac{1}{\Gamma(\beta)}\int_{0}^{+\infty}s^{\beta}e^{-s}u^{(k)}(x,t+s)\frac{ds}{s},\\
\end{eqnarray*}
 this implies, using Minkowski's integral inequality, 
\begin{eqnarray*}
\|\frac{\partial^{k}}{\partial
t^{k}} P_{t}(\mathcal{J}_{\beta}f)\|_{p,\gamma}&\leq&\frac{1}{\Gamma(\beta)}\int_{0}^{+\infty}s^{\beta}e^{-s}\|u^{(k)}(\cdot,t+s)\|_{p,\gamma}\frac{ds}{s}\\
&=& \frac{1}{\Gamma(\beta)}\int_{0}^{t}s^{\beta}e^{-s}\|u^{(k)}(\cdot,t+s)\|_{p,\gamma}\frac{ds}{s}\\
&&\quad \quad \quad + \frac{1}{\Gamma(\beta)}\int_{t}^{\infty}s^{\beta}e^{-s}\|u^{(k)}(\cdot,t+s)\|_{p,\gamma}\frac{ds}{s}\\
&=&(I)+(II).
\end{eqnarray*}
Now, as $\beta>0$, using Lemma \ref{kdecay} (as  $t+s>t$) and since
$f\in B_{p,\infty}^{\alpha}(\gamma_{d})$,
\begin{eqnarray*}
(I) &\leq&\frac{1}{\Gamma(\beta)}\|\frac{\partial^{k}P_{t}f}{\partial
t^{k}}\|_{p,\gamma}\int_{0}^{t}s^{\beta}e^{-s}\frac{ds}{s} \leq \frac{1}{\Gamma(\beta)}\|\frac{\partial^{k}P_{t}f}{\partial
t^{k}}\|_{p,\gamma}\int_{0}^{t}s^{\beta-1}ds\\
&\leq&\frac{1}{\Gamma(\beta)}\frac{t^{\beta}}{\beta}A_{k}(f)t^{-k+\alpha}
=C_{\beta} A_{k}(f)t^{-k+\alpha+\beta} .
\end{eqnarray*}
On the other hand, as $k>\alpha+\beta$  using Lemma
\ref{kdecay} as $t+s>s$, and since $f\in
B_{p,\infty}^{\alpha}(\gamma_{d})$
\begin{eqnarray*}
(II) &\leq&\frac{1}{\Gamma(\beta)}\int_{t}^{\infty}s^{\beta}e^{-s}\|\frac{\partial^{k}P_{s}f}{\partial
s^{k}}\|_{p,\gamma}\frac{ds}{s}
\leq \frac{A_{k}(f)}{\Gamma(\beta)}\int_{t}^{\infty}s^{\beta}e^{-s} s^{-k+\alpha}\frac{ds}{s}\\
&\leq& \frac{A_{k}(f)}{\Gamma(\beta)}\int_{t}^{\infty}s^{-k+\alpha+\beta-1}ds
=\frac{A_{k}(f)}{\Gamma(\beta)}\frac{t^{-k+\alpha+\beta}}{k-(\alpha+\beta)} = C_{k,\alpha,\beta}A_{k}(f) t^{-k+\alpha+\beta}.
\end{eqnarray*}
Therefore
\begin{eqnarray*}
\|\frac{\partial^{k}}{\partial
t^{k}} P_{t}({\mathcal{J}}_{\beta}f)\|_{p,\gamma}&\leq&C A_{k}(f) t^{-k+\alpha+\beta},
\end{eqnarray*}
 then
$\mathcal{J}_{\beta}f\in B_{p,\infty}^{\alpha+\beta}(\gamma_{d})$
and $A_{k}(\mathcal{J}_{\beta}f)\leq CA_{k}(f)$. Thus,
\begin{eqnarray*}
\|{\mathcal{J}}_{\beta}f\|_{B_{p,\infty}^{\alpha+\beta}}&=&\|{\mathcal{J}}_{\beta}f\|_{p,\gamma}+A_{k}({\mathcal{J}}_{\beta}f)\\
&\leq&\|f\|_{p,\gamma}+C A_{k}(f) \leq C\|f\|_{B_{p,\infty}^{\alpha}}. 
\end{eqnarray*}
\ep

In what follows we will need Hardy's inequalities, so for completeness we will write then here, see \cite{se70} page 272,
\begin{equation}\label{hardy1}
\int_{0}^{+\infty}\big(\int_{0}^{x}f(y)dy\big)^{p} x^{-r-1} dx \leq \frac{p}{r}\int_{0}^{+\infty}(y f(y))^{p}y^{-r-1}dy,
\end{equation}
and
\begin{equation}\label{hardy2}
\int_{0}^{+\infty}\big(\int_{x}^{\infty}f(y)dy\big)^{p} x^{r-1}dx \leq \frac{p}{r}\int_{0}^{+\infty}(y f(y))^{p}y^{r-1}dy,
\end{equation}
 where $f\geq 0, p\geq 1$ and $r>0.$\\

Now, we will study now the boundedness of the (Riesz) fractional derivative $D^\beta$ on Besov-Lipschitz spaces. We will use the representation (\ref{e6}) of the fractional derivative
and Hardy's inequalities.

\begin{teo}\label{DerRiesz<1}
Let $0<\beta<\alpha<1$, $1\leq p<\infty$ and $1 \leq q \leq \infty$ then
 $D^{\beta}$ is
bounded from $B_{p,q}^{\alpha}(\gamma_{d})$ into
$B_{p,q}^{\alpha-\beta}(\gamma_{d})$.
\end{teo}

\dem

 Let $f\in B_{p,q}^{\alpha}(\gamma_{d})$, using Hardy's inequality (\ref{hardy1}), with $p=1$, and the Fundamental Theorem of Calculus,
 \begin{eqnarray}\label{est4}
\nonumber |D^{\beta}f(x)| &\leq&\displaystyle\frac{1}{c_{\beta}}\int_{0}^{+\infty}s^{-\beta-1}|P_{s}f(x)-f(x)|ds\\
&\leq&
\nonumber \displaystyle\frac{1}{c_{\beta}}\int_{0}^{+\infty}s^{-\beta-1}\int_{0}^{s}|\frac{\partial
}{\partial r}P_{r}f(x)|dr \,ds\\
&\leq&\displaystyle\frac{1}{c_{\beta}\beta} \int_{0}^{+\infty}r^{1-\beta}|\frac{\partial
}{\partial r}P_{r}f(x)|\frac{dr}{r}  .
\end{eqnarray}
Thus, using Minkowski's integral inequality
\begin{eqnarray}\label{est5}
\|D^{\beta}f \|_{p,\gamma}
&\leq&\displaystyle C_{\beta} \int_{0}^{+\infty}r^{1-\beta}\|\frac{\partial}{\partial r}P_{r}f\|_{p,\gamma}\frac{dr}{r}<\infty ,
\end{eqnarray}
since  $f\in B_{p,q}^{\alpha}(\gamma_{d})\subset
B_{p,1}^{\beta}(\gamma_{d}) $, $1 \leq q \leq \infty$ as $\alpha> \beta,$  i.e. $D_{\beta}f \in L^{p}(\gamma_{d})  .$ \\

Now, by analogous argument,
\begin{eqnarray*}
\frac{\partial}{\partial t}P_{t}(D^{\beta}f)(x)&=& \displaystyle\frac{1}{c_{\beta}}\int_{0}^{+\infty}s^{-\beta-1}[\frac{\partial}{\partial t}P_{t+s}f(x)- \frac{\partial}{\partial t}P_t f(x)]ds\\
&=& \displaystyle\frac{1}{c_{\beta}}\int_{0}^{+\infty}s^{-\beta-1} \int_t^{t+s} u^{(2)}(x,r) dr \,ds\\
\end{eqnarray*}
and again, by Minkowski's integral inequality
\begin{equation}\label{pest1}
\| \frac{\partial}{\partial t}P_{t}(D^{\beta}f)\|\leq \displaystyle\frac{1}{c_{\beta}}\int_{0}^{+\infty}s^{-\beta-1} \int_t^{t+s} \|u^{(2)}(\cdot,r)\|_p dr \,ds\\
\end{equation}

Then, if $ 1 \leq q < \infty$, by (\ref{pest1})
\begin{eqnarray*}
&& \int_{0}^{\infty}\big(t^{1-(\alpha-\beta)}\|\frac{\partial }{\partial
t}P_{t}(D_{\beta}f)\|_{p,\gamma}\big)^{q}\frac{dt}{t}  \quad \quad \quad  \quad \quad \quad\\
& & \quad \quad \quad \leq \displaystyle C_{\beta}
\int_{0}^{\infty}\big(t^{1-(\alpha-\beta)}\int_{0}^{+\infty}s^{-\beta-1} \int_t^{t+s} \|u^{(2)}(\cdot,r) \|_{p,\gamma}dr \,ds\big)^{q}\frac{dt}{t}\\
 &&  \quad \quad \quad =\displaystyle C_{\beta}
\int_{0}^{\infty}\big(t^{1-(\alpha-\beta)}\int_{0}^{t}s^{-\beta-1} \int_t^{t+s} \|u^{(2)}(\cdot,r) \|_p dr \,ds\big)^{q}\frac{dt}{t}\\
&& \quad \quad \quad  \quad \quad + \, \displaystyle C_{\beta}
\int_{0}^{\infty}\big(t^{1-(\alpha-\beta)}\int_{t}^{+\infty}s^{-\beta-1} \int_t^{t+s} \|u^{(2)}(\cdot,r)\|_p dr \,ds\big)^{q}\frac{dt}{t}\\
&&  \quad \quad \quad =  (I)+(II).
\end{eqnarray*}
Now, since $r>t$ using Lemma \ref{kdecay} and the fact that  $0<\beta<1$,
\begin{eqnarray*}
(I) &\leq&C_{\beta} \int_{0}^{\infty}\big(t^{1-(\alpha-\beta)}\int_{0}^{t}s^{-\beta} ds\, \|u^{(2)}(\cdot,r)\|_{p,\gamma} \big)^{q}\frac{dt}{t}\\
&=& C_{\beta,q}\int_{0}^{\infty}\big(t^{2-\alpha}\|\frac{\partial^2}{\partial r^2}P_r f \|_{p,\gamma}\big)^{q}\frac{dt}{t}.
\end{eqnarray*}
On the other hand, as $r>t$ using Hardy's inequality (\ref{hardy2}), since $(1-\alpha)q>0$, we get
\begin{eqnarray*}
(II) &\leq& C_{\beta} \int_{0}^{\infty}t^{(1-(\alpha-\beta))q}\big(\int_{t}^{+\infty}s^{-\beta-1} ds \int_t^{\infty} \|u^{(2)}(\cdot,r) \|_{p,\gamma} dr \big)^{q}\frac{dt}{t}\\
&=& C'_{\beta} \int_{0}^{\infty}t^{(1-\alpha)q}\big( \int_t^{\infty} \|u^{(2)}(\cdot,r) \|_{p,\gamma} dr \, \big)^{q}\frac{dt}{t}\\
&\leq& \frac{C'_{\beta}}{(1-\alpha)} \int_{0}^{\infty}\big(r^{2-\alpha} \|\frac{\partial^2}{\partial r^2}P_r f \|_{p,\gamma} \big)^{q}\frac{dr}{r}.\\
\end{eqnarray*}
Thus,
\begin{eqnarray*}
\big(\int_{0}^{\infty}\big(t^{1-\alpha+\beta}\|\frac{\partial}{\partial t}P_{t}
D_{\beta}f\|_{p,\gamma}\big)^{q}\frac{dt}{t}\big)^{1/q}
&\leq&C \big(\int_{0}^{\infty}\big(t^{2-\alpha}\|\frac{\partial^{2}
}{\partial
t^{2}}P_{t}f\|_{p,\gamma}\big)^{q}\frac{dt}{t}\big)^{1/q}<\infty,
\end{eqnarray*}
as $f\in B_{p,q}^{\alpha}(\gamma_{d}).$ Then, $D_{\beta}f\in B_{p,q}^{\alpha-\beta}(\gamma_{d})$ and
\begin{eqnarray*}
\|D_{\beta}f\|_{B_{p,q}^{\alpha-\beta}}&=&\|D_{\beta}f\|_{p,\gamma}+\big(\int_{0}^{\infty}\big(t^{1-\alpha+\beta}\|\frac{\partial}{\partial t}P_{t}D_{\beta}f\|_{p,\gamma}\big)^{q}\frac{dt}{t}\big)^{1/q}\\
&\leq&C_{1} \|f\|_{B_{p,q}^{\alpha}}+C_{2} \big(\int_{0}^{\infty}\big(t^{2-\alpha}\|\frac{\partial^{2}
}{\partial t^{2}}P_{t}f\|_{p,\gamma}\big)^{q}\frac{dt}{t}\big)^{1/q}\\
&\leq&C \|f\|_{B_{p,q}^{\alpha}}  .
\end{eqnarray*}
Therefore $D_{\beta}f:B_{p,q}^{\alpha}\rightarrow
B_{p,q}^{\alpha-\beta}$ is bounded. \\

Now if $q=\infty$, inequality (\ref{pest1}) can be written as
\begin{eqnarray*}
\|\frac{\partial }{\partial t}P_{t}(D_{\beta}f)\|_{p,\gamma}
&\leq&\displaystyle\frac{1}{c_{\beta}}\int_{0}^{+\infty}s^{-\beta-1} \int_t^{t+s} \|\frac{\partial^2}{\partial r^2}P_r f \|_p dr \,ds\\
&=&\displaystyle\frac{1}{c_{\beta}}\int_{0}^{t}s^{-\beta-1} \int_t^{t+s} \|\frac{\partial^2}{\partial r^2}P_r f \|_p dr \,ds\\
&& \quad +\displaystyle\frac{1}{c_{\beta}}\int_{t}^{+\infty}s^{-\beta-1} \int_t^{t+s} \|\frac{\partial^2}{\partial r^2}P_r f \|_p dr \,ds=(I)+(II).
\end{eqnarray*}
Now,  by Lemma \ref{kdecay}, since $r>t$ 
\begin{eqnarray*}
(I) &\leq&\displaystyle\frac{1}{c_{\beta}}\int_0^t s^{-\beta}  \|\frac{\partial^2}{\partial t^2}P_t f \|_p  \,ds =C_\beta \|\frac{\partial^{2}
}{\partial t^{2}}P_{t}f\|_{p,\gamma} t^{1-\beta}\\
&\leq&C_\beta A(f) t^{-2+\alpha}t^{1-\beta}=C_\beta  A(f) t^{-1+\alpha-\beta},
\end{eqnarray*}
and by Lemma \ref{kdecay}, since $r>t$, and the fact that $f \in B_{p,\infty}^\alpha,$
\begin{eqnarray*}
(II)&\leq&\displaystyle\frac{1}{c_{\beta}}\int_{t}^{+\infty}s^{-\beta-1} \int_t^{\infty} \|\frac{\partial^2}{\partial r^2}P_r f \|_p dr \,ds\\
&\leq&C_\beta t ^{-\beta}\int_{t}^{\infty}\|\frac{\partial^{2}
}{\partial r^{2}}P_{r}f\|_{p,\gamma}dr
\leq C_\beta A(f) t^{-\beta} \int_{t}^{\infty}r^{-2+\alpha}dr\\
&=& C_{\alpha,\beta} A(f) t^{-1+\alpha-\beta} .
\end{eqnarray*}
Thus,
$$
\|\frac{\partial }{\partial t}P_{t}(D_{\beta}f)\|_{p,\gamma} \leq C A(f) t^{-1+\alpha-\beta},\quad t>0.$$

i.e., $D_{\beta}f \in B_{p,\infty}^{\alpha-\beta}(\gamma_{d})$ then
$A(D_{\beta}f)\leq C A(f)$, and 
\begin{eqnarray*}
\|D_{\beta}f\|_{B_{p,\infty}^{\alpha-\beta}}&=&\|D_{\beta}f\|_{p,\gamma}+A(D_{\beta}f)\\
&\leq&C_{1} \|f\|_{B_{p,\infty}^{\alpha}}+C_{2} A(f) \leq C \|f\|_{B_{p,\infty}^{\alpha}}  .
\end{eqnarray*}
Therefore $D_{\beta}:B_{p,\infty}^{\alpha}\rightarrow
B_{p,\infty}^{\alpha-\beta}$ is bounded. \ep \\

Now we will study now the boundedness of the Bessel fractional derivative on Besov-Lipschitz spaces, for $0<\beta<\alpha<1$ 
\begin{teo}\label{DerBessel<1}
Let $0<\beta<\alpha<1$, $1\leq p<\infty$ and $1 \leq q \leq \infty$ then
 ${\mathcal D}_{\beta}$ is
bounded from $B_{p,q}^{\alpha}(\gamma_{d})$ into
$B_{p,q}^{\alpha-\beta}(\gamma_{d})$.
\end{teo}

\dem

Let $f\in L^{p}(\gamma_{d})$, using the Fundamental Theorem of Calculus  we can write,
 \begin{eqnarray*}
\nonumber |\mathcal{D}_{\beta}f(x)| &\leq&\displaystyle\frac{1}{c_{\beta}}\int_{0}^{+\infty}s^{-\beta-1}|e^{-s}P_{s}f(x)-f(x)|ds\\
&\leq&\displaystyle\frac{1}{c_{\beta}}\int_{0}^{+\infty}s^{-\beta-1}e^{-s}|P_{s}f(x)-f(x)|ds
+\displaystyle\frac{1}{c_{\beta}}\int_{0}^{+\infty}s^{-\beta-1}|e^{-s}-1||f(x)|ds\\
&\leq&\displaystyle\frac{1}{c_{\beta}}\int_{0}^{+\infty}s^{-\beta-1}|\int_{0}^{s}\frac{\partial
}{\partial r}P_{r}f(x)dr|\,ds
+\displaystyle\frac{1}{c_{\beta}}\int_{0}^{+\infty}s^{-\beta-1}|e^{-s}-1||f(x)|ds\\
&\leq&\displaystyle\frac{1}{c_{\beta}}\int_{0}^{+\infty}s^{-\beta-1}\int_{0}^{s}|\frac{\partial
}{\partial r}P_{r}f(x)|dr\,ds
+\displaystyle\frac{1}{c_{\beta}}\int_{0}^{+\infty}s^{-\beta-1}|e^{-s}-1||f(x)|ds\\
&=&\displaystyle\frac{1}{c_{\beta}}\int_{0}^{+\infty}s^{-\beta-1}\int_{0}^{s}|\frac{\partial
}{\partial r}P_{r}f(x)|dr\,ds
+|f(x)|\displaystyle\frac{1}{c_{\beta}}\int_{0}^{+\infty}s^{-\beta-1}|-\int_{0}^{s}e^{-r}dr\,|ds\\
&=&\displaystyle\frac{1}{c_{\beta}}\int_{0}^{+\infty}s^{-\beta-1}\int_{0}^{s}|\frac{\partial
}{\partial r}P_{r}f(x)|dr\,ds
+\displaystyle\frac{1}{c_{\beta}} |f(x)|\int_{0}^{+\infty}s^{-\beta-1}\int_{0}^{s}e^{-r}dr\,ds.
\end{eqnarray*}

Now, using Hardy's inequality  (\ref{hardy1}), with $p=1$, in both integrals,  we have
 \begin{eqnarray*}
\nonumber |\mathcal{D}_{\beta}f(x)|
&\leq&\displaystyle\frac{1}{c_{\beta}}\int_{0}^{+\infty}s^{-\beta-1}\int_{0}^{s}|\frac{\partial
}{\partial r}P_{r}f(x)|dr\,ds
+\displaystyle\frac{1}{c_{\beta}}|f(x)|\int_{0}^{+\infty}s^{-\beta-1}\int_{0}^{s}e^{-r}dr \,ds\\
&\leq&\displaystyle\frac{1}{\beta c_{\beta}}\int_{0}^{+\infty}r|\frac{\partial
}{\partial r}P_{r}f(x)|r^{-\beta-1}dr
+\displaystyle\frac{1}{\beta c_{\beta}}|f(x)|\int_{0}^{+\infty}re^{-r}r^{-\beta-1}dr\\
&=&\displaystyle\frac{1}{\beta c_{\beta}}\int_{0}^{+\infty}r^{1-\beta}|\frac{\partial
}{\partial r}P_{r}f(x)|\frac{dr}{r}
+\displaystyle\frac{1}{\beta c_{\beta}}|f(x)|\int_{0}^{+\infty}r^{(1-\beta)-1}e^{-r}dr\\
&=& \displaystyle\frac{1}{\beta c_{\beta}}\int_{0}^{+\infty}r^{1-\beta}|\frac{\partial
}{\partial r}P_{r}f(x)|\frac{dr}{r}
+\displaystyle\frac{1}{\beta c_{\beta}}\Gamma(1-\beta)|f(x)|.
\end{eqnarray*}
Therefore, using the Minkowski's integral inequality
$$
\|\mathcal{D}_{\beta}f\|_p \leq
\displaystyle\frac{1}{\beta c_{\beta}}\int_{0}^{+\infty}r^{1-\beta}\|\frac{\partial
}{\partial r}P_{r}f\|_p\frac{dr}{r}
+\displaystyle\frac{1}{\beta c_{\beta}}\Gamma(1-\beta)\|f\|_p < C_1\|f\|_{B_{p,q}^{\alpha}}< \infty,
$$
since  $f\in B_{p,q}^{\alpha}(\gamma_{d})\subset
B_{p,1}^{\beta}(\gamma_{d}) $, $1 \leq q \leq \infty$ as $\alpha> \beta,$  i.e. $D_{\beta}f \in L^{p}(\gamma_{d}).$ \\

On the other hand, using the Fundamental Theorem of Calculus and using Hardy's inequality  (\ref{hardy1}), with $p=1$, in the second integral  we have,
 \begin{eqnarray*}
|\frac{\partial}{\partial t} P_{t}(\mathcal{D}_{\beta}f )(x)|&\leq&\frac 1{c_\beta}\int_0^{\infty}s^{-\beta-1} |e^{-s}\frac{\partial}{\partial t} P_{t+s}f(x)  -\frac{\partial}{\partial t} P_tf(x)| ds\\
&\leq&\frac 1{c_\beta}\int_0^{\infty}s^{-\beta-1}e^{-s}| \frac{\partial}{\partial t} P_{t+s}f(x)  -\frac{\partial}{\partial t} P_tf(x)| ds\\
&& \quad \quad +\frac 1{c_\beta}\int_0^{\infty}s^{-\beta-1}|e^{-s}  -1| |\frac{\partial}{\partial t} P_{t}f(x)| ds\\
&\leq&\frac 1{c_\beta}\int_0^{\infty}s^{-\beta-1} \int_t^{t+s} 
|\frac{\partial^2}{\partial r^2}P_{r} f(x)| dr \,ds\\
&& \quad  \quad + \frac 1{c_\beta} |\frac{\partial}{\partial t} P_{t}f(x)| \int_0^{\infty}s^{-\beta-1}\int_{0}^{s}e^{-r}dr\,ds,\\
&\leq&\frac 1{c_\beta}\int_0^{\infty}s^{-\beta-1} \int_t^{t+s} 
|\frac{\partial^2}{\partial r^2}P_{r} f(x)| dr \,ds\\
&& \quad \quad +\frac {1}{\beta c_\beta}|\frac{\partial}{\partial t} P_{t}f(x)|\int_0^\infty r^{(1-\beta) -1} e^{-r} dr\\
&=& \frac 1{c_\beta}\int_0^{\infty}s^{-\beta-1} \int_t^{t+s} 
|\frac{\partial^2}{\partial r^2}P_{r} f(x)| dr \,ds+\frac {\Gamma(1-\beta)}{\beta c_\beta}|\frac{\partial}{\partial t} P_{t}f(x)|.
\end{eqnarray*}
Therefore, by Minkowski's integral inequality 
\begin{equation}\label{pest2}
\|\frac{\partial}{\partial t} P_{t}(\mathcal{D}_{\beta}f )\|_{p,\gamma} \leq \frac 1{c_\beta}\int_0^{\infty}s^{-\beta-1} \int_t^{t+s} 
\|\frac{\partial^2}{\partial r^2}P_{r} f\|_{p,\gamma} dr \,ds+\frac {\Gamma(1-\beta)}{\beta c_\beta}\|\frac{\partial}{\partial t} P_{t}f\|_{p,\gamma}.
\end{equation}

Then, if  $1 \leq q < \infty$,  by (\ref{pest2})  Minkowski's integral inequality, we get
\begin{eqnarray*}
&&\big(\int_{0}^{\infty}\big(t^{1-(\alpha-\beta)}\|\frac{\partial }{\partial t}P_{t}\mathcal{D}_{\beta}f\|_{p,\gamma}\big)^{q}\frac{dt}{t}\big)^{1/q}\\
&&\quad \quad \quad\leq\frac {1}{ c_\beta}\big(\int_{0}^{\infty}\big(t^{1-(\alpha-\beta)} \int_0^{\infty}s^{-\beta-1} \int_t^{t+s} 
\|\frac{\partial^2}{\partial r^2}P_{r} f\|_{p,\gamma} dr \,ds \big)^{q}\frac{dt}{t}\big)^{1/q}\\
&&\quad \quad \quad \quad \quad + \frac {\Gamma(1-\beta)}{\beta c_\beta}\big(\int_{0}^{\infty}\big(t^{1-(\alpha-\beta)} \|\frac{\partial}{\partial t} P_{t}f\|_{p,\gamma}\big)^{q}\frac{dt}{t}\big)^{1/q}\\
&&\quad \quad \quad = (I) +(II).
\end{eqnarray*}
Now, the first term is the same as the one considered in the second part of the proof of Theorem \ref{DerRiesz<1}, thus by the same argument
\begin{eqnarray*}
(I)  &\leq&C_\beta\big(\int_{0}^{\infty}\big(t^{2-\alpha}\|\frac{\partial^{2}
}{\partial t^{2}}P_{t}f\|_{p,\gamma}\big)^{q}\frac{dt}{t}\big)^{1/q}< \|f\|_{B_{p,q}^{\alpha}}<\infty,
\end{eqnarray*}
since $f\in B_{p,q}^{\alpha}(\gamma_{d})$, and for the second term trivially
\begin{eqnarray*}
(II) &\leq&C\|f\|_{B_{p,q}^{\alpha-\beta}}\leq C\|f\|_{B_{p,q}^{\alpha}}
\end{eqnarray*}
since  $\alpha>\alpha- \beta$ and the inclusion relation, Proposition \ref{incluBesov}. \\

Thus  if  $1 \leq q < \infty$, 
$$\big(\int_{0}^{\infty}\big(t^{1-(\alpha-\beta)}\|\frac{\partial }{\partial t}P_{t}(\mathcal{D}_{\beta}f)\|_{p,\gamma}\big)^{q}\frac{dt}{t}\big)^{1/q} \leq C_2 \|f\|_{B_{p,q}^{\alpha}}$$
i.e. $\mathcal{D}_{\beta}f\in B_{p,q}^{\alpha-\beta}(\gamma_{d})$ and moreover
\begin{eqnarray*}
\|\mathcal{D}_{\beta}f\|_{B_{p,q}^{\alpha-\beta}}&=&\|\mathcal{D}_{\beta}f\|_{p,\gamma}+\big(\int_{0}^{\infty}\big(t^{1-\alpha+\beta}\|\frac{\partial
}{\partial
t}P_{t}\mathcal{D}_{\beta}f\|_{p,\gamma}\big)^{q}\frac{dt}{t}\big)^{1/q}\\
&\leq&C_{1}\|f\|_{B_{p,q}^{\alpha}}+C_{2}\big(\int_{0}^{\infty}\big(t^{2-\alpha}\|\frac{\partial^{2}
}{\partial
t^{2}}P_{t}f\|_{p,\gamma}\big)^{q}\frac{dt}{t}\big)^{1/q}\\
&\leq&C\|f\|_{B_{p,q}^{\alpha}}.
\end{eqnarray*}


If $q=\infty$, using the same argument as in Theorem \ref{DerRiesz<1}, inequality (\ref{pest2}) can be written as

\begin{eqnarray*}
\|\frac{\partial}{\partial t} P_{t}\mathcal{D}_{\beta}f \|_{p,\gamma}&\leq& \frac 1{c_\beta}\int_0^{\infty}s^{-\beta-1} \int_t^{t+s} 
\|\frac{\partial^2}{\partial r^2}P_{r} f\|_{p,\gamma} dr \,ds+\frac {\Gamma(1-\beta)}{\beta c_\beta}\|\frac{\partial}{\partial t} P_{t}f \|_{p,\gamma}\\
&\leq&C_{\alpha,\beta}A(f)t^{-1+\alpha-\beta}+\frac {\Gamma(1-\beta)}{\beta c_\beta}A(f)t^{-1+\alpha-\beta}\\
&\leq&C_{\alpha,\beta}A(f)t^{-1+\alpha-\beta},\quad t>0
\end{eqnarray*}
i.e. $\mathcal{D}_{\beta}f\in B_{p,\infty}^{\alpha-\beta}(\gamma_{d})$ and $A(\mathcal{D}_{\beta}f)\leq C_{\alpha,\beta}A(f)$, thus
\begin{eqnarray*}
\|\mathcal{D}_{\beta}f\|_{B_{p,\infty}^{\alpha-\beta}}&=&\|\mathcal{D}_{\beta}f\|_{p,\gamma}+A(\mathcal{D}_{\beta}f)\\
&\leq&C_{1} \|f\|_{B_{p,\infty}^{\alpha}}+C_{2} A(f)\leq C \|f\|_{B_{p,\infty}^{\alpha}}.  
\end{eqnarray*}
\ep \\

Now  we will consider the general case for fractional derivatives, removing the condition that the indexes must be less than 1. We need to consider forward differences. Remember for a given function $f$,  the $k$-th order forward difference of $f$ starting at $t$ with increment $s$ is defined as,\\
 $$\Delta_{s}^{k}(f,t)=\displaystyle\sum_{j=0}^{k}\binom{k}{j}(-1)^{j}f(t+(k-j)s).$$
 The forward differences have the following properties (see Appendix in \cite{gatur})
 we will need the following technical result
 \begin{lemma}\label{forw-diff} For any positive integer $k$
 \begin{enumerate}
 \item[i)]$\Delta_{s}^{k}(f,t)=\Delta_{s}^{k-1}(\Delta_{s}(f,\cdot),t)=\Delta_{s}(\Delta_{s}^{k-1}(f,\cdot),t)$
 \item[ii)] $\Delta_{s}^{k}(f,t)=\displaystyle\int_{t}^{t+s}\int_{v_{1}}^{v_{1}+s}...
 \int_{v_{k-2}}^{v_{k-2}+s}\int_{v_{k-1}}^{v_{k-1}+s}f^{(k)}(v_{k})dv_{k}dv_{k-1}...dv_{2}dv_{1}$
For any positive integer $k$,
\begin{equation}\label{difder}
\frac{\partial }{\partial s}(\Delta_s^k (f,t) ) = k \,\Delta_s^{k-1} (f',t+s),
\end{equation}
and for any integer $j>0$,
\begin{equation}\label{difder2}
\frac{\partial^j }{\partial t^j}(\Delta_s^k (f,t) ) =\Delta_s^{k} (f^{(j)},t).
\end{equation}
 \end{enumerate}
 \end{lemma}
 
 Observe that, using the Binomial Theorem and the semigroup property of $\{P_t\}$, we have 
 \begin{eqnarray}\label{powerrep}
\nonumber ( P_t -I )^k f(x) &=& \sum_{j=0}^k {k \choose j} P_t^{k-j} (-I)^j f(x) = \sum_{j=0}^k {k \choose j} (-1)^jP_t^{k-j} f(x)\\
\nonumber &=&\sum_{j=0}^k {k \choose j} (-1)^jP_{(k-j)t} f(x) =\sum_{j=0}^k {k \choose j} (-1)^j u(x,(k-j)t)\\
&=& \Delta_t^k (u(x, \cdot), 0),
\end{eqnarray}
where  as usual, $u(x,t) = P_t f(x)$.

Additionally we will need in what follows the following result,
 \begin{lemma}\label{lpforw-diff}
Let $f\in L^{p}(\gamma_{d}), \, 1\leq p<\infty$ and $k,n\in\mathbb{N}$ then
$$\displaystyle\|\Delta_{s}^{k}(u^{(n)},t)\|_{p,\gamma_{d}}\leq
s^{k}\|u^{(k+n)}(\cdot,t)\|_{p,\gamma_{d}} $$
 \end{lemma}
 \dem
 From ii) of Lemma \ref{forw-diff}, we have
  $$\Delta_{s}^{k}(u^{(n)}(x,\cdot),t)=\displaystyle\int_{t}^{t+s}\int_{v_{1}}^{v_{1}+s}...\int_{v_{k-2}}^{v_{k-2}+s}
 \int_{v_{k-1}}^{v_{k-1}+s}u^{(k+n)}(x,v_{k})dv_{k}dv_{k-1}...dv_{2}dv_{1},$$
 then, using Minkowski's integral inequality $k$-times and Lemma \ref{kdecay},
 \begin{eqnarray*}
\|\Delta_{s}^{k}(u^{(n)},t)\|_{p,\gamma_{d}}&\leq&\displaystyle\int_{t}^{t+s}\int_{v_{1}}^{v_{1}+s}...
\int_{v_{k-2}}^{v_{k-2}+s}\int_{v_{k-1}}^{v_{k-1}+s}\|u^{(k+n)}(\cdot,v_{k})\|_{p,\gamma_{d}}dv_{k}dv_{k-1}...dv_{2}dv_{1}\\
&\leq&s^{k}\|u^{(k+n)}(\cdot,t)\|_{p,\gamma_{d}}=s^{k}\|\frac{\partial^{k+n}}{\partial
 t^{k+n}}u(\cdot,t)\|_{p,\gamma_{d}}. 
 \end{eqnarray*}
\ep

Let us start with the case of the Riesz derivative,
\begin{teo}\label{DerRiesz>1}
Let $0<\beta<\alpha$, $1\leq p<\infty$ and $1\leq q\leq\infty$ then

 $D^{\beta}$ is
bounded from $B_{p,q}^{\alpha}(\gamma_{d})$ into
$B_{p,q}^{\alpha-\beta}(\gamma_{d})$.
\end{teo}
\dem

Let $f\in B_{p,q}^{\alpha}(\gamma_{d})$, using (\ref{powerrep}), Hardy's inequality (\ref{hardy1}), $p=1$, the Fundamental Theorem of Calculus and iii) of Lemma \ref{forw-diff}, we get
 \begin{eqnarray*}
 \nonumber |D^{\beta}f(x)|&\leq&\displaystyle\frac{1}{c_{\beta}}\int_{0}^{+\infty}s^{-\beta-1}|\Delta_{s}^{k}(u(x,\cdot),0)|ds\\
&\leq&\displaystyle\frac{1}{c_{\beta}}\int_{0}^{+\infty}s^{-\beta-1}\int_{0}^{s}|\frac{\partial}{\partial r}\Delta_{r}^{k}(u(x,\cdot),0)|dr\,ds\\
&\leq&\displaystyle\frac{1}{\beta c_{\beta}}\int_{0}^{+\infty}r^{-\beta}|\frac{\partial}{\partial r}\Delta_{r}^{k}(u(x,\cdot),0)|dr\\
&=&\displaystyle\frac{k}{\beta c_{\beta}}\int_{0}^{+\infty}r^{-\beta}|\Delta_{r}^{k-1}(u'(x,\cdot),r)|dr.
\end{eqnarray*}
Now, using Minkowski's integral inequality and Lemma \ref{lpforw-diff}
 \begin{eqnarray*}
 \|D_{\beta}f\|_{p,\gamma}&\leq&\displaystyle\frac{k}{\beta c_{\beta}}\int_{0}^{+\infty}
 r^{-\beta}\|\Delta_{r}^{k-1}(u',r)\|_{p,\gamma}dr\\
&\leq&\displaystyle\frac{k}{\beta c_{\beta}}\int_{0}^{+\infty}
 r^{k-\beta}\|\frac{\partial^{k}}{\partial
 r^{k}}P_{r}f\|_{p,\gamma}\frac{dr}{r}<\infty,
 \end{eqnarray*}
since $f\in B_{p,q}^{\alpha}(\gamma_{d})\subset
B_{p,1}^{\beta}(\gamma_{d})$, as $\alpha>\beta$. Therefore, $D_{\beta}f\in L^{p}(\gamma_{d})$.\\

On the other hand, 
\begin{eqnarray*}
P_t [( P_s -I )^k f(x) ]  &=& P_t (\Delta_s^k (u(x, \cdot), 0)) = P_t ( \sum_{j=0}^k {k \choose j} (-1)^jP_{(k-j)s} f(x) )\\
&=&    \sum_{j=0}^k {k \choose j} (-1)^jP_{t+(k-j)s} f(x) = \Delta_s^k (u(x, \cdot), t).
\end{eqnarray*}
Thus, if $n$ be the smaller integer greater than $\alpha$, i.e. $ n-1 \leq \alpha  < n$, then by Lemma \ref{forw-diff} iv),
 \begin{eqnarray*}
 \frac{\partial^{n}}{\partial t^{n}}P_{t}(D_{\beta}f)(x)&=&\frac{1}{c_{\beta}}\int_{0}^{+\infty}s^{-\beta-1}\frac{\partial^{n}}{\partial t^{n}}(\Delta_{s}^{k}(u(x,\cdot),t)\\
 &=&\frac{1}{c_{\beta}}\int_{0}^{+\infty}s^{-\beta-1}\Delta_{s}^{k}(u^{(n)}(x,\cdot),t)ds.
 \end{eqnarray*}
 and therefore, by Minkowski's integral inequality
 \begin{equation}\label{pest3}
\| \frac{\partial^{n}}{\partial t^{n}}P_{t}(D_{\beta}f)\|_{p,\gamma} \leq \frac{1}{c_{\beta}}\int_{0}^{+\infty}s^{-\beta-1}\| \Delta_{s}^{k}(u^{(n)},t)\|_{p,\gamma} ds.
\end{equation}

Now if $1 \leq q< \infty$, by  (\ref{pest3}), 
  \begin{eqnarray*}
&&  \big(\int_{0}^{\infty}\big(t^{n-(\alpha-\beta)}\|\frac{\partial^{n}}{\partial t^{n}}P_{t}(D_{\beta}f)\|_{p,\gamma}\big)^{q}\frac{dt}{t}\big)^{1/q}\\
&& \quad \quad \quad \leq \frac{1}{c_{\beta}} \big(\int_{0}^{\infty}\big(t^{n-(\alpha-\beta)}\int_{0}^{+\infty}s^{-\beta-1}\|\Delta_{s}^{k}(u^{(n)},t)
\|_{p,\gamma}ds\big)^{q}\frac{dt}{t}\big)^{1/q}\\
&& \quad \quad \quad \leq \frac{1}{c_{\beta}}\big(\int_{0}^{\infty}\big(t^{n-(\alpha-\beta)}\int_{0}^{t}s^{-\beta-1}\|\Delta_{s}^{k}(u^{(n)},t)
\|_{p,\gamma}ds\big)^{q}\frac{dt}{t}\big)^{1/q}\\
&& \quad \quad \quad  \quad \quad + \frac{1}{c_{\beta}} \big(\int_{0}^{\infty}\big(t^{n-(\alpha-\beta)}\int_{t}^{+\infty}s^{-\beta-1}\|\Delta_{s}^{k}(u^{(n)},t)
\|_{p,\gamma}ds\big)^{q}\frac{dt}{t}\big)^{1/q}\\
&& \quad \quad \quad=(I)+(II).\\
\end{eqnarray*}

Then, by Lemma \ref{lpforw-diff},
 \begin{eqnarray*}
 (I) &\leq& \frac{1}{c_{\beta}}\big(\int_{0}^{\infty}\big(t^{n-(\alpha-\beta)}\|\frac{\partial^{n+k}}{\partial t^{n+k}}P_{t}f\|_{p,\gamma}\int_{0}^{t}s^{k-\beta-1}ds\big)^{q}\frac{dt}{t}\big)^{1/q}\\
&=&\frac{1}{c_{\beta}(k-\beta)}\big(\int_{0}^{\infty}\big(t^{n+k-\alpha}\|u^{(n+k)}(\cdot,t)
\|_{p,\gamma}\big)^{q}\frac{dt}{t}\big)^{1/q}<\infty,
\end{eqnarray*}
since $f\in B_{p,q}^{\alpha}(\gamma_{d})$,
and by Lemma \ref{kdecay}
\begin{eqnarray*}
 (II)
&\leq&\frac{1}{c_{\beta}}\big(\int_{0}^{\infty}\big(t^{n-(\alpha-\beta)}\int_{t}^{+\infty}s^{-\beta-1}\big(\sum_{j=0}^{k}\binom{k}{j}\|u^{(n)}(\cdot,t+(k-j)s)\|_{p,\gamma}\big)ds\big)^{q}\frac{dt}{t}\big)^{1/q}\\
&\leq&\frac{1}{c_{\beta}}\big(\int_{0}^{\infty}\big(t^{n-(\alpha-\beta)}\int_{t}^{+\infty}s^{-\beta-1}\big(\sum_{j=0}^{k}\binom{k}{j}\|u^{(n)}(\cdot,t)
\|_{p,\gamma}\big)ds\big)^{q}\frac{dt}{t}\big)^{1/q}\\
&=&\frac{2^{k}}{c_{\beta}}\big(\int_{0}^{\infty}\big(t^{n-(\alpha-\beta)}\|\frac{\partial^{n}}{\partial t^{n}}P_{t}f\|_{p,\gamma}\int_{t}^{+\infty}s^{-\beta-1}ds\big)^{q}\frac{dt}{t}\big)^{1/q}\\
&=&\frac{2^{k}}{c_{\beta}\beta}\big(\int_{0}^{\infty}\big(t^{n-\alpha}\|\frac{\partial^{n}}{\partial t^{n}}P_{t}f\|_{p,\gamma}\big)^{q}\frac{dt}{t}\big)^{1/q}<\infty,
\end{eqnarray*}
since $f\in B_{p,q}^{\alpha}(\gamma_{d})$. Therefore,  if $ 1 \leq q < \infty$, $D_{\beta}f\in B_{p,q}^{\alpha-\beta}(\gamma_{d})$, and moreover,
\begin{eqnarray*}
\|D_{\beta}f\|_{B_{p,q}^{\alpha-\beta}}&=&\|D_{\beta}f\|_{p,\gamma}+\big(\int_{0}^{\infty}\big(t^{n-\alpha+\beta}\|\frac{\partial^{n}}{\partial t^{n}}P_{t}(D_{\beta}f)\|_{p,\gamma}\big)^{q}\frac{dt}{t}\big)^{1/q}\\
&\leq&C_{1}\|f\|_{B_{p,q}^{\alpha}}+C_{2}\|f\|_{B_{p,q}^{\alpha}}\leq C\|f\|_{B_{p,q}^{\alpha}}
\end{eqnarray*}
Thus, $D_{\beta}f:B_{p,q}^{\alpha}\rightarrow
B_{p,q}^{\alpha-\beta}$ is bounded.\\

If $q =\infty$, inequality (\ref{pest3}) can be written as 
\begin{eqnarray*}
\| \frac{\partial^{n}}{\partial t^{n}}P_{t}(D_{\beta}f)\|_{p,\gamma} &\leq& \frac{1}{c_{\beta}}\int_{0}^{t}s^{-\beta-1}\| \Delta_{s}^{k}(u^{(n)},t)\|_{p,\gamma} ds\\
 && \quad \quad \quad +\frac{1}{c_{\beta}}\int_{t}^{+\infty}s^{-\beta-1}\| \Delta_{s}^{k}(u^{(n)},t)\|_{p,\gamma} ds\\
 &=& (I) + (II)
\end{eqnarray*}
and then as $f \in B^\alpha_{p, \infty}$,  by Lemma \ref{lpforw-diff},
\begin{eqnarray*}
(I) &\leq& \frac{1}{c_{\beta}}\int_{0}^{t}s^{-\beta-1} s^k \|u^{(n+k)} \|_{p,\gamma} ds = C_\beta \|\frac{\partial^{n+k}}{\partial t^{n+k}}P_{t}f\|_{p,\gamma} t^{k-\beta} \\
&\leq& C_\beta A(f) t^{-n-k+\alpha} t^{k-\beta} = C_\beta A(f) t^{-n+\alpha-\beta}, 
\end{eqnarray*}
and  as above, by Lemma \ref{kdecay},
\begin{eqnarray*}
(II) &\leq& \frac{1}{c_{\beta}} \int_{t}^{+\infty}s^{-\beta-1}\big(\sum_{j=0}^{k}\binom{k}{j}\|u^{(n)}(\cdot,t+(k-j)s)\|_{p,\gamma}\big)ds\\
&\leq& C_\beta \int_{t}^{+\infty}s^{-\beta-1}\big(\sum_{j=0}^{k}\binom{k}{j}\|u^{(n)}(\cdot,t)\|_{p,\gamma}) ds= C_\beta t^{-\beta} \|\frac{\partial^{n}}{\partial t^{n}}P_{t}f\|_{p,\gamma}\\
&\leq& C_\beta A(f) t^{-n+\alpha} t^{-\beta} = C_\beta A(f) t^{-n+\alpha-\beta}. 
\end{eqnarray*}

 \ep \\

There is an alternative proof of the fact that  $D_{\beta}f\in L^{p}(\gamma_{d})$ without using Hardy's inequality following the same scheme as in the proof of i) Theorem 3.5 in  \cite{gatur}, using the inclusion $B^\alpha_{p,q} \subset B^{\beta+\epsilon}_{p, \infty}$ with $\beta + \epsilon <k$.\\

\begin{teo}\label{DerBessel>1}
Let $0<\beta<\alpha$, $1\leq p<\infty$ and $1\leq q\leq\infty$ then

 $\mathcal{D}_{\beta}$ is
bounded from $B_{p,q}^{\alpha}(\gamma_{d})$ into
$B_{p,q}^{\alpha-\beta}(\gamma_{d})$.
\end{teo}
\dem

Let $f\in B_{p,q}^{\alpha}(\gamma_{d})$, and set $v(x,t)=e^{-t} u(x,t)$ then using the Hardy's inequality (\ref{hardy1}), the Fundamental Theorem of Calculus and iii) of Lemma \ref{forw-diff},
 \begin{eqnarray*}
 \nonumber |\mathcal{D}_{\beta}f(x)|&\leq&\displaystyle\frac{1}{c_{\beta}}\int_{0}^{+\infty}s^{-\beta-1}|\Delta_{s}^{k}(v(x,\cdot),0)|ds\\
&\leq&\displaystyle\frac{1}{c_{\beta}}\int_{0}^{+\infty}s^{-\beta-1}\int_{0}^{s}|\frac{\partial}{\partial r}\Delta_{r}^{k}(v(x,\cdot),0)|dr\,ds\\
&\leq&\displaystyle\frac{k}{\beta c_{\beta}}\int_{0}^{+\infty}r^{-\beta}|\Delta_{r}^{k-1}(v'(x,\cdot),r)|dr
\end{eqnarray*}
and this implies by Minkowski's integral inequality 
$$\|\mathcal{D}_{\beta}f\|_{p,\gamma_{d}}\leq\displaystyle\frac{k}{\beta c_{\beta}}
\int_{0}^{+\infty}r^{-\beta}\|\Delta_{r}^{k-1}(v',r)\|_{p,\gamma}dr.$$
Now, using Lemma \ref{forw-diff}

$\|\Delta_{r}^{k-1}(v',r)\|_{p,\gamma}\leq\displaystyle\int_{r}^{2r}\int_{v_{1}}^{v_{1}+r}...
 \int_{v_{k-2}}^{v_{k-2}+r}\|v^{(k)}(\cdot,v_{k-1})\|_{p,\gamma}dv_{k-1}dv_{k-2}...dv_{2}dv_{1}$
and by Leibnitz's differentiation rule for the product
 \begin{eqnarray*}
\|v^{(k)}(\cdot,v_{k-1})\|_{p,\gamma}&=&\displaystyle\|\sum_{j=0}^{k}\binom{k}{j}(e^{-v_{k-1}})^{(j)}u^{(k-j)}(\cdot,v_{k-1})\|_{p,\gamma_{d}}\\
&\leq&\displaystyle\sum_{j=0}^{k}\binom{k}{j}e^{-v_{k-1}}\|u^{(k-j)}(\cdot,v_{k-1})\|_{p,\gamma}.
  \end{eqnarray*}
  Then
  \begin{eqnarray*}
&& \|\Delta_{r}^{k-1}(v',r)\|_{p,\gamma}\\
&& \quad \quad \quad \leq\displaystyle\sum_{j=0}^{k}\binom{k}{j}\int_{r}^{2r}\int_{v_{1}}^{v_{1}+r}...
 \int_{v_{k-2}}^{v_{k-2}+r}e^{-v_{k-1}}\|u^{(k-j)}(\cdot,v_{k-1})\|_{p,\gamma}dv_{k-1}dv_{k-2}...dv_{2}dv_{1}\\
 &&\quad \quad \quad \leq \sum_{j=0}^{k}\binom{k}{j}r^{k-1}e^{-r}\|u^{(k-j)}(\cdot,r)\|_{p,\gamma}.
 \end{eqnarray*}
Therefore
\begin{eqnarray*}
 \nonumber \|\mathcal{D}_{\beta}f\|_{p,\gamma}&\leq&\displaystyle\frac{k}{\beta c_{\beta}}
\sum_{j=0}^{k}\binom{k}{j}\int_{0}^{+\infty}
 r^{k-\beta-1}e^{-r}\|u^{(k-j)}(\cdot,r)\|_{p,\gamma}dr\\
&=&\displaystyle\frac{k}{\beta c_{\beta}}\sum_{j=0}^{k-1}\binom{k}{j}\int_{0}^{+\infty}
 r^{(k-j)-(\beta-j)-1}e^{-r}\|\frac{\partial^{k-j}}{\partial r^{k-j}}P_{r}f\|_{p,\gamma}dr\\
&& \quad \quad \quad +\frac{k}{\beta c_{\beta}}\int_{0}^{+\infty}
 r^{k-\beta-1}e^{-r}\|P_{r}f\|_{p,\gamma}dr\\
 &\leq&\displaystyle\frac{k}{\beta c_{\beta}}\sum_{j=0}^{k-1}\binom{k}{j}\int_{0}^{+\infty}
 r^{(k-j)-(\beta-j)-1}\|\frac{\partial^{k-j}}{\partial r^{k-j}}P_{r}f\|_{p,\gamma}dr\\
&& \quad \quad \quad + \frac{k}{\beta c_{\beta}}\int_{0}^{+\infty}
 r^{k-\beta-1}e^{-r}\|f\|_{p,\gamma}dr
\end{eqnarray*}
 Thus
\begin{eqnarray*}
 \nonumber \|\mathcal{D}_{\beta}f\|_{p,\gamma}&\leq&\displaystyle\frac{k}{\beta c_{\beta}}\sum_{j=0}^{k-1}\binom{k}{j}\int_{0}^{+\infty}
 r^{k-j-(\beta-j)}\|\frac{\partial^{k-j}}{\partial r^{k-j}}P_{r}f\|_{p,\,\gamma}\frac{dr}{r}\\
&& \quad \quad \quad\quad \quad \quad \quad \quad \quad + \frac{k \Gamma(k-\beta)}{\beta c_{\beta}}\|f\|_{p,\gamma}<\infty,
\end{eqnarray*}
 since $f\in B_{p,q}^{\alpha}(\gamma_{d})\subset B_{p,1}^{\beta-j}(\gamma_d) $ as $\alpha>\beta>\beta-j\geq 0$, for $j\in\{0,...,k-1\}$, then $\mathcal{D}_{\beta}f\in L^{p}(\gamma_d)$.\\

On the other hand,
 $$P_{t}(e^{-s}P_{s}-I)^{k}f(x)=\displaystyle\sum_{j=0}^{k}\binom{k}{j}(-1)^{j}e^{-s(k-j)}u(x,t+(k-j)s).$$
Let $n$ be the smaller integer greater than $\alpha$, i.e. $ n-1 \leq \alpha  < n$, we have
 \begin{eqnarray*}
 \frac{\partial^{n}}{\partial t^{n}}P_{t}(\mathcal{D}_{\beta}f)(x)&=&\frac{1}{c_{\beta}}
 \int_{0}^{+\infty}s^{-\beta-1}\sum_{j=0}^{k}\binom{k}{j}(-1)^{j}e^{-s(k-j)}u^{(n)}(x,t+(k-j)s)ds\\
 &=&\frac{e^t}{c_{\beta}}
 \int_{0}^{+\infty}s^{-\beta-1}\sum_{j=0}^{k}\binom{k}{j}(-1)^{j}e^{-(t+s(k-j))}u^{(n)}(x,t+(k-j)s)ds\\
&=&\frac{e^{t}}{c_{\beta}}
 \int_{0}^{+\infty}s^{-\beta-1}\Delta_{s}^{k}(w(x,\cdot),t) ds,
  \end{eqnarray*}
where $w(x,t)=\displaystyle e^{-t}u^{(n)}(x,t)$. Now using the Fundamental Theorem of Calculus,
 \begin{eqnarray*}
 \frac{\partial^{n}}{\partial t^{n}}P_{t}(\mathcal{D}_{\beta}f)(x)&=&\frac{e^{t}}{c_{\beta}}
 \int_{0}^{+\infty}s^{-\beta-1}\Delta_{s}^{k}(w(x,\cdot),t)ds\\
&=&\frac{e^{t}}{c_{\beta}}
  \int_{0}^{+\infty}s^{-\beta-1}\int_{0}^{s}\frac{\partial}{\partial
  r}\Delta_{r}^{k}(w(x,\cdot),t)dr\,ds.
  \end{eqnarray*}
Then, using Hardy's inequality (\ref{hardy1}), and iii) of Lemma \ref{forw-diff},
 \begin{eqnarray*}
 |\frac{\partial^{n}}{\partial t^{n}}P_{t}\mathcal({D}_{\beta}f)(x)|&\leq&\frac{e^{t}}{c_{\beta}}
  \int_{0}^{+\infty}s^{-\beta-1}\int_{0}^{s}|\frac{\partial}{\partial
  r}\Delta_{r}^{k}(w(x,\cdot),t)|drds\\
&\leq&\frac{e^{t}}{c_{\beta}\beta}
  \int_{0}^{+\infty}r|\frac{\partial}{\partial
  r}\Delta_{r}^{k}(w(x,\cdot),t)|r^{-\beta-1}dr\\
  &=&\frac{ke^{t}}{c_{\beta}\beta}
  \int_{0}^{+\infty}r^{-\beta}|\Delta_{r}^{k-1}(w'(x,\cdot),t+r)|dr
  \end{eqnarray*}
  and by Minkowski's integral inequality we get
 \begin{eqnarray*}
 \|\frac{\partial^{n}}{\partial t^{n}}P_{t}(\mathcal{D}_{\beta}f)\|_{p,\gamma}&\leq&\frac{ke^{t}}{\beta c_{\beta}}
  \int_{0}^{+\infty}r^{-\beta}\|\Delta_{r}^{k-1}(w',t+r)\|_{p,\gamma}dr.
  \end{eqnarray*}
  Now, by analogous argument as above, Lemma \ref{forw-diff} and Leibnitz's pro-duct rule give us
 \begin{eqnarray*}
\|\Delta_{r}^{k-1}(w',t+r)\|_{p,\gamma}&\leq&\sum_{j=0}^{k}\binom{k}{j}r^{k-1}e^{-(t+r)}\|u^{(k+n-j)}(\cdot,t+r)\|_{p,\gamma},
\end{eqnarray*}
and this implies that
 \begin{eqnarray*}\label{pest4}
 \|\frac{\partial^{n}}{\partial t^{n}}P_{t}(\mathcal{D}_{\beta}f)\|_{p,\gamma}&\leq&e^{t}\frac{k}{c_{\beta}\beta}
  \int_{0}^{+\infty}r^{-\beta}\big(\sum_{j=0}^{k}\binom{k}{j}r^{k-1}e^{-(t+r)}\|u^{(k+n-j)}(\cdot,t+r)\|_{p,\gamma}\big)
dr\\
 &=&\frac{k}{c_{\beta}\beta}
  \sum_{j=0}^{k}\binom{k}{j}\int_{0}^{+\infty}r^{k-\beta-1}e^{-r}\|u^{(k+n-j)}(\cdot,t+r)\|_{p,\gamma}
dr.
\end{eqnarray*}
Thus 
\begin{equation}
 \|\frac{\partial^{n}}{\partial t^{n}}P_{t}(\mathcal{D}_{\beta}f)\|_{p,\gamma} \leq 
 \frac{k}{c_{\beta}\beta}
  \sum_{j=0}^{k}\binom{k}{j}\int_{0}^{+\infty}r^{k-\beta-1}e^{-r}\|u^{(k+n-j)}(\cdot,t+r)\|_{p,\gamma}\, dr.
\end{equation}

Now if $1 \leq q< \infty$, using (\ref{pest4}) we have,
 \begin{eqnarray*}
 &&\big(\int_{0}^{\infty}\big(t^{n-(\alpha-\beta)}\|\frac{\partial^{n}}{\partial t^{n}}P_{t}(\mathcal{D}_{\beta}f)\|_{p,\gamma}\big)^{q}\frac{dt}{t}\big)^{1/q}\\
 &&\quad \quad \leq\frac{k}{c_{\beta}\beta}
  \sum_{j=0}^{k}\binom{k}{j} \big(\int_{0}^{\infty}\big(t^{n-(\alpha-\beta)}\int_{0}^{+\infty}r^{k-\beta-1}e^{-r}\|u^{(k+n-j)}(\cdot,t+r)\|_{p,\gamma}
dr\big)^{q}\frac{dt}{t}\big)^{1/q}.
\end{eqnarray*}
For each $1\leq j\leq k$,
$0<\alpha-\beta+k-j\leq\alpha$ and  by Lemma  \ref{kdecay}
\begin{eqnarray*}
&&
(\int_{0}^{\infty}\big(t^{n-(\alpha-\beta)}\int_{0}^{\infty}r^{k-\beta-1}e^{-r}\|u^{(k+n-j)}(\cdot,t+r)\|_{p,\gamma}
dr\big)^{q}\frac{dt}{t}\big)^{1/q}\\
&& \quad \quad \quad  \leq
(\int_{0}^{\infty}\big(t^{n-(\alpha-\beta)}\|u^{(n+k-j)}(\cdot,t)\|_{p,\gamma}\int_{0}^{+\infty}r^{k-\beta-1}e^{-r}
dr\big)^{q}\frac{dt}{t}\big)^{1/q}\\
&& \quad \quad \quad  = \Gamma(k-\beta)
(\int_{0}^{\infty}\big(t^{n+(k-j)-(\alpha-\beta+k-j)}\|u^{(n+k-j)}(\cdot,t)\|_{p,\gamma}\big)^{q}\frac{dt}{t}\big)^{1/q}<\infty,
\end{eqnarray*}
as  $f\in B_{p,q}^{\alpha}(\gamma_{d})\subset
B_{p,q}^{\alpha-\beta+(k-j)}(\gamma_{d})$ for any $0 \leq j \leq k$.\\

Now, for the case $j=0$,
 \begin{eqnarray*}
&&\big(\int_{0}^{\infty}\big(t^{n-(\alpha-\beta)}\int_{0}^{+\infty}r^{k-\beta-1}e^{-r}\|u^{(n+k)}(\cdot,t+r)\|_{p,\gamma}
dr\big)^{q}\frac{dt}{t}\big)^{1/q}\\
&&\quad \quad \quad \leq\big(\int_{0}^{\infty}\big(t^{n-(\alpha-\beta)}\int_{0}^{t}r^{k-\beta-1}e^{-r}\|u^{(n+k)}(\cdot,t+r)\|_{p,\gamma}
dr\big)^{q}\frac{dt}{t}\big)^{1/q}\\
&&\quad \quad \quad \quad \quad +\big(\int_{0}^{\infty}\big(t^{n-(\alpha-\beta)}\int_{t}^{+\infty}r^{k-\beta-1}e^{-r}\|u^{(n+k)}(\cdot,t+r)\|_{p,\gamma}
dr\big)^{q}\frac{dt}{t}\big)^{1/q}\\
&&\quad \quad \quad =(I)+(II).
\end{eqnarray*}
Using Lemma \ref{kdecay}, and $k>\beta$, 
 \begin{eqnarray*}
(I) &\leq&\big(\int_{0}^{\infty}\big(t^{n-(\alpha-\beta)}\int_{0}^{t}r^{k-\beta-1}\|u^{(n+k)}(\cdot,t)\|_{p,\gamma}
dr\big)^{q}\frac{dt}{t}\big)^{1/q}\\
&=&\big(\int_{0}^{\infty}\big(t^{n-(\alpha-\beta)}\|u^{(n+k)}(\cdot,t)\|_{p,\gamma}\int_{0}^{t}r^{k-\beta-1}dr\big)^{q}\frac{dt}{t}\big)^{1/q}\\
&=&\frac{1}{k-\beta}\big(\int_{0}^{\infty}\big(t^{n+k-\alpha}\|u^{(n+k)}(\cdot,t)\|_{p,\gamma}\big)^{q}\frac{dt}{t}\big)^{1/q}<\infty,
\end{eqnarray*}
since $f\in B_{p,q}^{\alpha}(\gamma_{d})$ and $n+k>\alpha$ and for the second term, using Lemma \ref{kdecay} and Hardy's inequality  (\ref{hardy2})
 \begin{eqnarray*}
(II) &\leq&\big(\int_{0}^{\infty}\big(t^{n-(\alpha-\beta)}\int_{t}^{+\infty}r^{k-\beta-1}\|u^{(n+k)}(\cdot,r)\|_{p,\gamma}
dr\big)^{q}\frac{dt}{t}\big)^{1/q}\\
&\leq&\frac{1}{n-(\alpha-\beta)}\big(\int_{0}^{\infty}\big(r^{n+k-\alpha}\|u^{(n+k)}(\cdot,r)\|_{p,\gamma}\big)^{q} \frac{dr}{r}\big)^{1/q}<\infty,
\end{eqnarray*}
since $f\in B_{p,q}^{\alpha}(\gamma_{d})$.\\

Therefore $\mathcal{D}_{\beta}f\in
B_{p,q}^{\alpha-\beta}(\gamma_{d})$. Moreover
 \begin{eqnarray*}
\|\mathcal{D}_{\beta}f\|_{B_{p,q}^{\alpha-\beta}}&=&\|\mathcal{D}_{\beta}f\|_{p,\gamma}+
\big(\int_{0}^{\infty}\big(t^{n-(\alpha-\beta)}\|\frac{\partial^{n}}{\partial t^{n}}P_{t}\mathcal{D}_{\beta}f\|_{p,\gamma}\big)^{q}\frac{dt}{t}\big)^{1/q}\\
 &\leq&C_{1}\|f\|_{p,\gamma}+\frac{k}{c_{\beta}\beta}\sum_{j=0}^{k}\binom{k}{j} C_{2}\big(\int_{0}^{\infty}\big(r^{n-\alpha}\|\frac{\partial^{n}}{\partial r^{n}}P_{r}f\|_{p,\gamma}\big)^{q}\frac{dr}{r}\big)^{1/q}\\
&\leq&C\|f\|_{B_{p,q}^{\alpha}}
\end{eqnarray*}

Finally, if $q= \infty$, from the inequality (\ref{pest4})
\begin{eqnarray*}
\|\frac{\partial^{n}}{\partial t^{n}}P_{t}(\mathcal{D}_{\beta}f)\|_{p,\gamma}\leq \frac{k}{c_{\beta}\beta}
  \sum_{j=0}^{k}\binom{k}{j}\int_{0}^{+\infty}r^{k-\beta-1}e^{-r}\|u^{(k+n-j)}(\cdot,t+r)\|_{p,\gamma}
dr,
\end{eqnarray*}
and then, the argument is essentially similar to the previous case, as we did  in the last part of the proof of Theorem \ref{DerRiesz>1}. 
 
 \ep \\


{\bf Observation}
Let us observe that  if instead of considering the {\em Ornstein-Uhlenbeck operator} (\ref{OUop}) and the {\em Poisson-Hermite semigroup} (\ref{PoissonH}) we consider the {\em  Laguerre differential  operator} in $\mathbb R^d_{+}$,
\begin{equation}
\mathcal{L}^{\alpha} = \sum^{d}_{i=1} \bigg[ x_i \frac{\partial^2}{\partial x^2_i}
+ (\alpha_i + 1 - x_i )  \frac{\partial}{\partial x_i} \bigg],
\end{equation}
and the corresponding {\em Poisson-Laguerre semigroup}, or if we consider the {\em Jacobi differential  operator} in $(-1,1)^d$,
\begin{equation}
\mathcal{L}^{\alpha,\beta} = - \sum^{d}_{i=1} \bigg[ (1-x_i^2)\frac{\partial^2}{\partial x_i^2}
+  (\beta_i -\alpha_i-\left(\alpha_i +\beta_i +2\right)x_i) \frac{\partial}{\partial x_i} \bigg],
\end{equation}
and the corresponding {\em Poisson-Jacobi semigroup} (for details we refer to \cite{ur2}), the arguments are completely analogous.  That is to say, we can defined in analogous manner {\em Laguerre-Besov-Lipschitz spaces, } and {\em Jacobi-Besov-Lipschitz spaces} then prove that the corresponding notions of Fractional Integrals and Fractional Derivatives  behave similarly.
In order to see this it is more convenient to use the representation (\ref{PoissonH}) of $P_t$ in terms of the one-sided stable measure $\mu^{(1/2)}_t(ds)$, see \cite{piur02}.



\begin{thebibliography}{99}

\bibitem{ButzBer}
 Butzer, P. L. \& Berens, H. \emph{Semi-groups of operators and approximation.} Die Grundkehren der mathematischen Wissenschaften, Ban145. Springer-Verlag. New York, 1967.
\bibitem{fm91}
  Frazier M, Jawerth B, Weiss G. \emph{Littlewood Paley Theory and the Study of Functions Spaces.} CBMS-Conference Lecture Notes 79.
 Amer. Math. Soc. Providence RI , 1991.
 \bibitem{forscotur}
Forzani, L., Scotto, R, and Urbina, W{.} \emph{Riesz and Bessel
Potentials, the $g_k$ functions and an Area function, for  the
Gaussian measure $\gamma_d$.} Revista de la Uni\'on Matem\'atica
 Argentina (UMA), vol 42 (2000), no.1,17--37.
\bibitem{gs96}
 Gatto A. E, Segovia C, V\'{a}gi S. \emph{On Fractional Differentiation and Integration on Spaces of Homogeneous Type.} Rev. Mat. Iberoamericana,  {\bf 12} (1996), 111--145.
  \bibitem{gatur}
 Gatto A. E, Urbina W.(2009) \emph{On Gaussian Lipschitz Spaces and the Boundedness of Fractional Integrals and Fractional Derivatives on them}. Preprint. arXiv:0911.3962v2.
\bibitem{lu}
L\'{o}pez I{.} and  Urbina, W{.} \emph{Fractional Differentiation
for the Gaussian Measure and Applications. } Bull. Sciences Math, 2004,
{\bf 128}, 587--603.
\bibitem{me3}
 Meyer, P. A. {\em Transformations de Riesz pour les lois Gaussiennes.}
 Lectures Notes in Math  1059 (1984) Springer-Verlag 179-193.
\bibitem{ebner} Pineda, E. \emph{T\'opicos en An\'alisis Arm\'onico Gaussiano: Comportamiento en la frontera y  Espacios de funciones para la medida Gaussiana}.
 Doctoral Thesis, Facultad de Ciencias, UCV, Caracas.(2009)
 \bibitem {piur02}
 Pineda, E. and Urbina, W. \emph{Some results on Gaussian Besov-Lipschitz and Gaussian Triebel-Lizorkin spaces.} Journal of Approximation Theory Vol 161, \# 2 (2009),  529-564.
 \bibitem{sam}
  Samko S., Kilbas, A \& Marichev, O.  \emph{Fractional integrals and derivatives: theory and applications.} Gordon and Breach Science Publishers, Philadelphia, 1992.
\bibitem{se70}
  Stein E. \emph{Singular integrals and differentiability properties of functions.} Princeton Univ. Press. Princeton, New Jersey, 1970.
     \bibitem{trie}
 Triebel, H  \emph{Interpolation theory, function spaces differential operators.} Noth Holland, 1978.
 \bibitem{trie1}
 Triebel, H  \emph{Theory of function spaces.} Birkh\"auser Verlag, Basel, 1983.
   \bibitem{trie2}
 Triebel, H  \emph{Theory of function spaces II.} Birkh\"auser Verlag, Basel, 1992.
  \bibitem{ur2}
 Urbina, W.  {\em Operators Semigroups associated to Classical Orthogonal Polynomials  and  Functional Inequalities.} Lecture Notes of the French Mathematical Society (SMF). 2008.
\bibitem{uw98}
  Urbina W. \emph{An\'{a}lisis Arm\'{o}nico Gaussiano: una visi\'{o}n panor\'{a}mica}. Trabajo de Ascenso, Facultad de Ciencias, UCV, Caracas, 1998.
  \bibitem{wat}
 Watanabe, S. {\em Lecture on Stochastic Differential Equations and 
Malliavin Calculus.} Tata Institute. Springer Verlag (1984).

\end{thebibliography}
\end{document}